\documentclass[final]{siamltex}

\usepackage{amsmath}
\usepackage{amsfonts}
\usepackage{graphics}
\usepackage{amssymb}
\usepackage{psfrag}
\usepackage{amsmath}
\usepackage{graphicx}
\usepackage{multirow}
\usepackage{bm}
\usepackage{cite}
\usepackage{bm}
\usepackage{marginnote}
\usepackage{mathtools}
\usepackage{footnote}
\usepackage{epic}

\usepackage{amsmath,amsfonts}
\usepackage{subfig}
\usepackage{setspace}
\usepackage{algpseudocode}
\usepackage{cite}
\usepackage{amsmath}
\usepackage{amssymb}
\usepackage{graphicx}
\usepackage{xcolor}
\usepackage{epstopdf}
\usepackage{bm}
\usepackage{tikz}
\usepackage{tabularx}
\usepackage{array}

\newcommand{\f}{\frac}
\newcommand{\ds}{\displaystyle}

\def\be{\begin{equation}}
\def\ee{\end{equation}}

\def\f{\frac}
\def\E{\mathbb{E}}
\def\ben{\begin{eqnarray}}
\def\een{\end{eqnarray}}

\begin{document}

\title{An efficient meshfree implicit filter for nonlinear filtering problems\thanks{This material is based upon work supported in part by the U.S.~Air Force of Scientific Research under grant numbers 1854-V521-12 and FA9550-11-1-0149; by the U.S.~Department of Energy, Office of Science, Office of Advanced Scientific Computing Research, Applied Mathematics program under contract and award numbers ERKJ259, ERKJE45, and DE-SC0010678; and by the Laboratory Directed Research and Development program at the Oak Ridge National Laboratory, which is operated by UT-Battelle, LLC., for the U.S.~Department of Energy under Contract DE-AC05-00OR22725.}}

\pagestyle{myheadings}
\markboth{F.~Bao, Y.~Cao, C.~G.~Webster and G.~Zhang}
         {A meshfree implicit filter for nonlinear filtering problems}

\author {Feng Bao\thanks{Department of Computational and Applied Mathematics, Oak Ridge National Laboratory, Oak Ridge, TN 37831 ({\tt baof@ornl.gov}, {\tt webstercg@ornl.gov}, {\tt zhangg@ornl.gov}).}
\and Yanzhao Cao\thanks{Department of Mathematics and Statistics, Auburn University, Auburn, Alabama, 36849 ({\tt yzc0009@auburn.edu}).}
\and Clayton G.~Webster\footnotemark[2]
\and Guannan Zhang\footnotemark[2]
}

\maketitle
\begin{abstract}
In this paper, we propose a meshfree approximation method for the implicit filter developed in \cite{Bao-implicit}, which is a novel numerical algorithm for nonlinear filtering problems. The implicit filter approximates conditional distributions in the optimal filter over a deterministic state space grid and is developed from samples of the current state obtained by solving the state equation implicitly.
The purpose of the meshfree approximation is to improve the efficiency of the implicit filter in moderately high-dimensional problems. The construction of the algorithm includes generation of random state space points and a meshfree interpolation method. Numerical experiments show the effectiveness and efficiency of our algorithm.
\end{abstract}

\begin{keywords}
Nonlinear filtering, implicit algorithm, meshfree approximation, Shepard's method
\end{keywords}

\section{Introduction}
Nonlinear filters are important tools for dynamical data assimilation with applications in a variety of research areas, including  
biology \cite{Baker2013, Little-Jones2013}, mathematical finance \cite{Bensoussan2009, Elliott2013},
signal processing \cite{Hairer2011, Little-Jones2013, Stannat2011}, image processing \cite{Singh2013}, 
and multi-target tracking \cite{Kim2013,Yang2012}. To put it succinctly, nonlinear filtering is an extension of the Bayesian framework to the estimation and prediction of nonlinear stochastic dynamics.
In this effort, we consider the following nonlinear filtering model
\begin{equation}\label{State:Cont}
\left\{
\begin{aligned}
\f{d X_t}{dt} &= f(t, X_t; W_t),  \quad \text{(state)}\\
 Y_t &= g(t, X_t) + V_t, \;\;\; \text{(observation)}
\end{aligned}\right.
\end{equation}
where $f$ and $g$ are two nonlinear functions, $\{X_t \in \mathbb{R}^d,t \geq 0 \}$ and $\{{Y}_t \in \mathbb{R}^q, t \geq 0\}$ are the stochastic state and observation processes, respectively, $\{W_t \in \mathbb{R}^{r}, t \ge 0\}$ is a random vector representing the uncertainty in $X_t$, and $\{V_t \in \mathbb{R}^{s}, t\ge 0\}$ denotes the random measurement error in $Y_t$. In the discrete setting, the nonlinear filtering model in \eqref{State:Cont} takes the form
%
\begin{equation}\label{eq:State}
\left\{
\begin{aligned}
{X}_{k} &= f_{k}({X}_{k-1}, w_{k-1}),  \quad \text{(state)}\\
{Y}_{k} &= g_{k}({X}_{k}) + v_{k}, \;\;\qquad \text{(observation)}
\end{aligned}\right.
\end{equation}
where $\{w_k\}_{k \in \mathbb{N}^{+}} \in \mathbb{R}^{r}$ and $\{v_k\}_{k \in \mathbb{N}^{+} } \in \mathbb{R}^{s}$ are mutually independent white noises. 
Let ${Y}_{1:k} := \sigma\{{Y}_1, {Y}_2, \cdots, {Y}_k \}$ denote the $\sigma$ filed generated by the observational data up to the step $k$. 
The goal of nonlinear filtering is to find the posterior probability density function (PDF) of the state ${X}_k$, given the observation data $Y_{1:k}$, so as to compute the quantity of interest (QoI), given by
$$ \mathbb{E}[\Phi({X}_k) | Y_{1:k} ] = \inf\left\{ \mathbb{E}[ | \Phi({X}_k) - Z |^2 ]; Z \in \mathcal{Z}_k \right\},$$
where $\Phi(\cdot)$ is a test function, and $\mathcal{Z}_k$ denotes the space of all $Y_{1:k}$-measurable and square integrable random variables.

Tremendous efforts have been made to solve nonlinear filtering problems in the last few decades.
Two of the well-known Bayesian filters are extended Kalman filters (EKFs) \cite{EKF-tracking, Dun-EKF, EKF, Julier-EKF, Kulikov-EKF}, and particle filters \cite{particle-filter-resample, Crisan-PF, Crisan-Xiong-PF, particle-filter}. 
The key ingredient of the EKFs is the linearization of both $f$ and $g$ in \eqref{State:Cont}, so that the standard Kalman filter can be applied directly. Thus, if the nonlinearity of the state and the observation systems is not severe, then the EKFs can provide efficient and reasonable inferences about the state, otherwise, the performance of the EKFs can be very poor. 
For particle filters, the central theme is to approximate the desired posterior PDF of the state by the empirical distribution of a set of adaptively selected random samples (referred to as ``particles''). The particle filter method is essentially a sequential Monte Carlo approach, which requires no assumption on the linearity of the underlying system. As such, with sufficiently large number of particles, it is capable of providing an accurate approximation of the posterior PDF for a highly nonlinear filtering problems. However, there are some fundamental issues concerning the efficiency and robustness of particle filters \cite{cd2002}. For example, since the empirical PDF is constructed based on particles with equal weights after resampling, the particle filter still needs a lot of samples in order to accurately approximate the target distribution. 

To overcome such a disadvantage, the authors proposed a new nonlinear filter named ``implicit filter" \cite{Bao-implicit} . This approach adopts the framework of  Bayesian filtering, which has two stages at each time step, i.e., prediction and update. At the prediction stage, we estimate the prior PDF $p(X_k | Y_{1:k-1} )$ of the future state $X_{k}$ given the current available observation information $Y_{1:k-1}$; at the update stage, we update the prior PDF by assimilating the newly received data $Y_{k}$ to obtain the estimate of the posterior PDF $p(X_k | Y_{1:k})$. The implicit filter is distinguished from the particle filters by the use of interpolatory approximations to the prior and posterior PDFs. Specifically, in the particle filter, $p(X_k | Y_{1:k-1} )$ is approximated by {\em explicitly} propagating the samples of the current state $X_{k-1}|Y_{1:k-1}$ through the nonlinear state equation ${X}_{k} = f_{k}({X}_{k-1}, w_{k-1})$, and constructing the empirical PDF of $X_{k}|Y_{1:k-1}$. In the implicit filter, the interpolation of $p(X_k | Y_{1:k-1} )$ requires its function values at a set of grid points of the future state $X_k$. Under the condition that $X_k = x \in \mathbb{R}^d$, we solve {\em implicitly} the state equation $x = f_{k}({X}_{k-1}, w_{k-1})$ given a set of Monte Carlo samples of $w_{k-1}$, so that the value of $p(X_k=x | Y_{1:k-1} )$, at the grid point of $x$, can be estimated by averaging the function values of $p(X_{k-1}|Y_{1:k-1})$ at all the solutions of the state equation. As an implicit scheme, the implicit filter has a stabilizing effect which provides more accurate numerical approximations to the solution of the nonlinear filtering problem than the particle filter method \cite{Bao-implicit}.

%

The main challenge of the implicit filter method is that the conditional PDF of the nonlinear filtering solution is estimated at grid points. As such the method suffers the so called  ``the curse of dimensionality'' when the dimension of the state variable is high. In addition,  the efficiency of the method may be significantly reduced when the domain of the PDF is unbounded. 
In this paper, we propose to construct a meshfree implicit filter algorithm to alleviate the aforementioned challenges. 
Motivated by the particle filter method, we first generate a set of random particles and  propagate these particles  through the system model and use these particles to replace the grid points in the state space. After that we generate other necessary points through the  Shepard's method  which constructs the interpolant by the weighted average of the values on state points \cite{Fasshauer2007}.  In order to prevent  particle degeneracy in the generation of random state points, we introduce a resample step in the particle propagation.  In addition  we choose  state points  according to the system state, which make them  adaptively located in the high probability region of the PDF of state. In this way, we solve the nonlinear filtering problem in a relatively small region in the state space at each time step and approximate the solution on a set of meshfree state points distributed adaptively to the desired PDF of the state. Furthermore, since we approximate the PDF as a function on each state point, instead of using state points themselves to describe the empirical distribution, the implicit filter algorithm requires much fewer points than the particle filter method to depict the PDF of the state.  

The rest of this paper is organized as follows. In \S \ref{BF}, we introduce the mathematical framework of the Bayesian optimal filter. In \S \ref{Algorithm}, we construct meshfree implicit algorithm. In \S \ref{sec:ex}, we demonstrate the efficiency and accuracy of our algorithm through numerical experiments. Finally, \S \ref{sec:con} contains conclusions and directions for the future research.


\section{Bayesian optimal filter}\label{BF}
For $m,n \in \mathbb{N}^{+}$, 
let $X_{m:n}$ and $Y_{m:n}$ denote the $\sigma$ fields generated by $\{X_m, X_{m+1}, \ldots, X_n\}$ and $\{Y_m, Y_{m+1}, \ldots, Y_n\}$, respectively. For $k = \mathbb{N}^{+}$, we use $x_k$ to represent a realization of the random variable $X_k$, and define
\[
p(x_{k}|\cdot) := p( X_k=x_k|\cdot)
\]
for notational simplicity. 
It is easy to see that the  dynamical model in \eqref{eq:State} is Markovian in the sense that 
$$ 
p(x_k | X_{1 : k-1}, Y_{1 : k-1}) = p(x_k | X_{k-1}).
$$
We also know that the measurements are conditionally independent given $x_k$, i.e.,
$$p(Y_k | X_{1 : k}, Y_{1 : {k-1}}) = p(Y_k | x_k). $$
The Bayesian optimal filter constructs the conditional distribution $p( x_k | Y_{1:k} )$ recursively in two stages: prediction stage  and update stage.

For $k=1,  2, \cdots$,   assume that  $p(x_{k-1} | Y_{1:k-1})$ is given.  In the prediction stage  $p(x_k | Y_{1:k-1})$ is evaluated  through   the Chapman-Kolmogorov formula: 
\begin{equation}
p(x_k | Y_{1:k-1}) = \int_{\mathbb{R}^{d}} p(x_k | x_{k-1}) p(x_{k-1} | Y_{1:k-1}) d x_{k-1}.  \label{prediction_bayes}
\end{equation}
In the update stage,   the prior PDF obtained in \eqref{prediction_bayes} is used to obtain the posterior PDF $p(x_k | Y_{1:k} )$   via the Bayes' formula: 
\begin{equation}
p(x_k | Y_{1:k} ) = \f{p(Y_k | x_k) p(x_k | Y_{1:k-1} )}{p(Y_k | Y_{1:k-1} )} = \f{p(Y_k | x_k) p(x_k | Y_{1:k-1} )}{\int_{\mathbb{R}^{d}} p(Y_k | x_k) p(x_k | Y_{1:k-1} )\, d x_{k}}.  \label{bayes}
\end{equation}

\section{The meshfree implicit filter}\label{Algorithm}
In this section, we construct the  meshfree implicit filter algorithm.  The algorithm is based the implicit filter algorithm on grid points \cite{Bao-implicit}. 
The implicit filter algorithm introduced in \cite{Bao-implicit} is developed from the general framework of the Bayesian optimal filter discussed above, in which the primary computational challenge is the numerical approximation of the term $  p(x_k | x_{k-1}) $ in \eqref{prediction_bayes}. 

\subsection{The prediction stage}\label{Prediction}
For $k=1, 2, \cdots$, the goal of this stage is to approximate the prior distribution $p(x_k | Y_{1:k-1})$ of the state $X_k$, given the posterior distribution $p(x_{k-1} | Y_{1:k-1} )$ of the state $X_{k-1}$. Due to the the fact that 
$$
p(x_k | x_{k-1}) = \mathbb{E}_{w}[p(x_k | x_{k-1} , w_{k-1})] = \int_{\mathbb{R}^{r}} p(x_k | x_{k-1}, w_{k-1}) \cdot p(w_{k-1}) d w_{k-1} ,
$$
the prior PDF $p(x_k | Y_{1:k-1} ) $ derived in identity \eqref{prediction_bayes} can be rewritten as
\begin{equation} \label{predic}
\begin{aligned}
p(x_k | Y_{1:k-1} ) = \int_{\mathbb{R}^{d}} \mathbb{E}_w[p(x_k | x_{k-1} , w_{k-1})] p(x_{k-1} | Y_{1:k-1} ) d x_{k-1},
\end{aligned}
\end{equation}
where $\mathbb{E}_w[\cdot]$ represents the expectation with respect to the white noise $w_{k-1}$, and the PDF $p(x_k | x_{k-1}, w_{k-1})$ is 
\begin{equation}\label{ppp}
p(x_k | x_{k-1}, w_{k-1}) = \left\{
\begin{aligned}
\infty, \;\;\; x_k = f_k(x_{k-1}, w_{k-1}),\\
0, \;\;\; x_k \neq  f_k(x_{k-1}, w_{k-1}),
\end{aligned}
\right.
\end{equation}
with $\int_{\mathbb{R}^d} p(x_k | x_{k-1}, w_{k-1}) d{x_{k}} = 1$ for any $x_{k-1} \in \mathbb{R}^d$ and $w_{k-1}\in \mathbb{R}^r$. The definition in \eqref{ppp}
can be viewed as a generalization of the Dirac delta function in the space $\mathbb{R}^d \times \mathbb{R}^d \times \mathbb{R}^r$, where the mass is located according to the state equation $x_k = f_k(x_{k-1}, w_{k-1})$.

Note that the estimation of \eqref{predic} requires an approximation to the expectation $ \mathbb{E}_w[p(x_k | x_{k-1} , w_{k-1})] $. To this end, we first draw $M$ independent samples $\{w_{k-1}^{j}\}_{j = 1}^M$ of the white noise $w_{k-1}$, and define an approximation to $p(x_k | x_{k-1}, w_{k-1})$ as
\begin{equation}\label{PPPI}
\begin{aligned}
\pi^{M} (x_k| x_{k-1}, w_{k-1}) & := \sum_{j=1}^M \delta_{w_{k-1}^{j}}(x_k| x_{k-1}, w_{k-1}), \\
\end{aligned}
\end{equation}
with 
\[
\delta_{w_{k-1}^{j}}(x_k| x_{k-1}, w_{k-1}) := \left\{
\begin{aligned}
\infty, &\;\;\; w_{k-1} = w_{k-1}^j \text{ and } x_k = f_k(x_{k-1}, w_{k-1}^j),\\
0, &\;\;\; \text{otherwise},\\
\end{aligned}
\right.
\]
which is essentially a restriction of $p(x_k| x_{k-1}, w_{k-1})$ in the subset $\{w_{k-1}^{j}\}_{j = 1}^M$. 
%
%
%
Therefore, the expectation $\E_w[p(x_k | x_{k-1} , w_{k-1})]$ in \eqref{predic} can be approximated by
\begin{equation}\label{Exp:empirical}
\begin{aligned}
\E_w[p(x_k | x_{k-1} ,w_{k-1})] & \approx \mathbb{E}_w\left[\pi^M(x_k | x_{k-1} , w_{k-1})\right], \\
& =  \sum_{j=1}^M \int_{\mathbb{R}^r} \delta_{w_{k-1}^{j}}(x_k| x_{k-1}, w_{k-1}) p(w_{k-1})dw_{k-1}.
\end{aligned}
\end{equation}

To construct an interpolation of $p( x_k | Y_{1:k-1})$, the next step is to approximate $p( x_k | Y_{1:k-1})$ at a point set 
$\mathcal{H}_k := \{x_k^{i}\}_{i=1}^N \subset \mathbb{R}^d$ with $N \in \mathbb{N}^+$.
By substituting $x_k = x_k^i$ into \eqref{predic}-\eqref{Exp:empirical}, we have
\begin{equation}
p( x_k^i | Y_{1:k-1} ) =  \int_{\mathbb{R}^{d}} \E_w\left[\pi^M(  x_k^i | x_{k-1} , w_{k-1})\right] p(x_{k-1} | Y_{1:k-1} ) d x_{k-1} + \mathcal{R}_{k | k-1}^i, \label{integ}
\end{equation}
where $\mathcal{R}_{k | k-1}^i := p( x_k^i | Y_{1:k-1} ) - \int_{\mathbb{R}^{d}} \E_w[\pi^M(  x_k^i | x_{k-1} , w_{k-1})] p(x_{k-1} | Y_{1:k-1} ) d x_{k-1}$ is the approximation error. Then, by further fixing $w_{k-1} = w_{k-1}^{j}$, the location of the mass of $\delta_{w_{k-1}^{j}}(x_k^i| x_{k-1}, w_{k-1}^j)$ in the space of $x_{k-1}$, denoted by $x^{i,j}_{k-1}$, can be obtained by
{\em implicitly} solving the state equation 
$$
f_k\left(x^{i,j}_{k-1}, w_{k-1}^{j}\right) =  x_k^i , \quad j = 1, \cdots, M,
$$
which is the reason we named the approach the implicit filter. 
%
Now substituting $x_{k}^i$ into \eqref{Exp:empirical}, and using the same sample set $\{w_{k-1}^{j}\}_{j=1}^M$ as above to approximate the integral on the right hand side of \eqref{Exp:empirical},
we obtain
\begin{equation}\label{bbb}
\begin{aligned}
\E_w\left[\pi^M\left(  x_k^i | x_{k-1}, w_{k-1}\right)\right] & =  \sum_{j=1}^M \left(\f{1}{M}\sum_{j'=1}^M \delta_{w_{k-1}^{j}}(x_k^i| x_{k-1}, w_{k-1}^{j'})\right)\\
& = \f{1}{M}\sum_{j=1}^M \delta_{w_{k-1}^{j}}\left(x_k^i| x_{k-1}, w_{k-1}^{j}\right),
\end{aligned}
\end{equation}
then replacing $\E_w[\pi^M(  x_k^i | x_{k-1} , w_{k-1})] $ in   \eqref{integ} with \eqref{bbb}, we have  
\begin{equation}\label{scheme:predict}
 \begin{aligned}
 p(  x_k^i | Y_{1:k-1} ) & =  \ds \int_{\mathbb{R}^{d}}  \left[\f{1}{M}\sum_{j=1}^M \delta_{w_{k-1}^{j}}\left(x_k^i| x_{k-1}, w_{k-1}^{j}\right)\right] p(x_{k-1} | Y_{1:k-1} ) d x_{k-1} + \mathcal{R}_{k | k-1}^i\\
   & = \f{1}{M}\sum^M_{j=1} p \left( x^{i,j}_{k-1} \Big| Y_{1:k-1} \right) + \mathcal{R}_{k | k-1}^i,
  \end{aligned}
\end{equation}
where $ p( x^{i, j}_{k-1} | Y_{1:k-1})$ is the value of $p \left( x_{k-1} | Y_{1:k-1} \right)$ at $x^{i, j}_{k-1}$.
Neglecting the error term $\mathcal{R}_{k | k-1}^i$ in \eqref{scheme:predict}, we obtain the following iterative numerical scheme for constructing an approximation, denoted by $\varrho(x_{k}^i | Y_{1:k-1} )$, of the prior PDF $p(  x_k^i | Y_{1:k-1} )$, i.e.,
\begin{equation}\label{scheme:a}
\varrho(x_{k}^i | Y_{1:k-1} ) =   \f{1}{M}\sum^M_{j=1} \varrho ( x^{i,j}_{k-1} | Y_{1:k-1} ).
\end{equation}

In our previous work \cite{Bao-implicit}, the subsets $\mathcal{H}_k$, for $k = 0, 1, \ldots$, were defined by a full tensor product mesh, denoted by
 \begin{equation}\label{mesh:a}
\mathcal{M} := \mathcal{M}^{(1)} \times \mathcal{M}^{(2)} \times \cdots \mathcal{M}^{(d)},
\end{equation}
on a $d$-dimensional hyper-cube $[a_1, b_1] \times \cdots \times [a_d, b_d]$, where $\mathcal{M}^{(m)}, m = 1, \dots, d$, is a uniform partition of
the interval $[a_m, b_m]$ with $N^{(m)}$ grid points. It is simple to implement but has several significant disadvantages. First, at each time step, one needs to approximate the prior PDF $p(  x_k | Y_{1:k-1} )$ at a total of $N^{(1)} \times \cdots \times N^{(d)}$ grid points which grows exponentially as the dimension $d$ increases. This is also known as ``the curse of dimensionality''. On the other hand, since the construction of $\mathcal{M}$ is not informed by the target PDF,  
the domain $[a_1, b_1] \times \cdots \times [a_d, b_d]$ needs to be defined sufficiently large, so as to capture the statistically significant region of the PDF. This may lead to a great waste of computation effort in the low probability region of $p(  x_k | Y_{1:k-1} )$.

To alleviate such disadvantages, we propose to develop a distribution-informed meshfree interpolation approach to efficiently approximate the prior PDF. 
%
%
The central idea of the generation of random points for the state variable is to build a set of points, denoted by  $\mathcal{H}_k$,   according to the state distribution. 
To begin with, we generate  $\mathcal{H}_0=\{\xi^{i}\}_{i = 1}^{N}$ of $N$ random samples from the initial PDF $p_0$ of the initial state: 
$$\mathcal{H}_0 := \{x_0^{i}\}_{i = 1}^N = \{\xi^{i}\}_{i=1}^N, \ \text{with} \ x_0^{i} = \xi^{i} .$$ 
If the initial PDF $p_0$ is close to the true state distribution, it's obvious that our random state points are more concentrated near the target state. For $k = 1, 2, \cdots, K$, we propagate points $\{ x_{k-1}^{i}\}_{i = 1}^N$ to $\{ x_{k}^{i}\}_{i = 1}^N$ through the state equation \eqref{eq:State}:
$$ x_{k}^{i} = f_{k-1}(x_{k-1}^{i}, \tilde{w}_{k-1}^{i}), \quad  i = 1, 2, \cdots, N, $$
where $\{\tilde{w}_{k-1}^{i}\}_{i = 1}^N$ are $N$ random samples according to the PDF of $w_{k-1}$.  Denote $\mathcal{H}_k := \{x_k^{i}\}_{i = 1}^N $ and approximate the conditional PDF $p(x_k | Y_{1:k-1} )$ on $\mathcal{H}_k$ with the scheme given by \eqref{scheme:a}.
In this way, the random points in $\mathcal{H}_k$ move according to the state model. 
As opposed to particle filter methods, which use the number of particles to represent empirical distributions and require a large number of particles to follow the state distribution, in the implicit filter method we provide an approximation of the value of the PDF at each state point. Therefore, much fewer points are needed to describe the state PDF and the random state points are not necessary to accurately follow the state distribution.

\subsection{The update stage}\label{Update}
By incorporating the new data $Y_{k}$, we update the prior PDF $p(x_k | Y_{1:k-1})$ at each grid point $x_{k}^i$, using the Bayesian formula, to obtain
\begin{equation}\label{scheme:update}
\begin{aligned}
p(  x_k^i | Y_{1:k} )  =& \ \f{1}{C_k} p(Y_k |  x_k^i ) p(x_k^i | Y_{1:k-1} )  \\
=& \ \f{1}{C_k}  p(Y_k |  x_k^i) \varrho(x_k^i | Y_{1:k-1} )  + \mathcal{R}_{k | k}^i,
\end{aligned}
\end{equation}
where $\varrho(x_k^i | Y_{1:k-1} ) $ is given in \eqref{scheme:a}, $C_k$ is the normalization factor, and $\mathcal{R}_{k | k}^i :=  \f{1}{C_k}  p(Y_k |  x_k^i)  \big( p(x_k^i | Y_{1:k-1} ) - \varrho(x_k^i | Y_{1:k-1} ) \big) $ is the approximation error. By neglecting the error term $\mathcal{R}_{k | k}^i$ in $\eqref{scheme:update}$, we obtain the following iterative numerical scheme for the update stage on $\mathcal{H}_k$, i.e.,
\begin{eqnarray}
\varrho(  x_k^i | Y_{1:k} ) &=& \f{1}{C_k} p(Y_k |  x_k^i) \varrho(x_k^i | Y_{1:k-1} )  \label{scheme:b},
\end{eqnarray}
where $\varrho(  x_k^i | Y_{1:k})$ is desired the approximation of the posterior PDF $p(  x_k^i | Y_{1:k} )$.

Next, we use interpolation methods to construct the approximation $\varrho(x_{k} | Y_{1:k} )$ of $p(x_{k} | Y_{1:k} )$ from values $\{\varrho( x_k^i | Y_{1:k} )\}_{x_k^i \in \mathcal{H}_k }$ via
\begin{equation}\label{Interpolation}
\varrho(x_{k} | Y_{1:k} ) = \sum_{x_k^i \in \mathcal{H}_k} \varrho(  x_k^i | Y_{1:k} ) \phi^{i}(x_k),
\end{equation}
where $\{\phi^{i}\}_{i=1}^{N_k}$ is the set of basis functions. 
Since the state points in $\mathcal{H}_k$ are generated randomly in the meshfree framework, standard polynomial interpolation \cite{Bao-implicit} is unstable due to the uncontrollable Lebesgue constant. Instead, we propose to use the Shepard's method \cite{Fasshauer2007}, which is an efficient meshfree interpolation technique, to construct the interpolant $\varrho(x_{k} | Y_{1:k} )$. The basic idea of the Shepard's method is to use the weighted average of  $\{\varrho(x_{k}^i | Y_{1:k} )\}_{x_{k}^i \in \mathcal{H}_{k}}$ in the interpolating approximation. Specifically, for a given point $x_{k} \in \mathbb{R}^{d}$, we re-order the points in $\mathcal{H}_{k}$ by the distances to $x_{k}$ to get a sequence $\{ x_{k}^{(l)} \}_{l = 1}^{N_k}$ such that
$$ \|x_{k} - x_{k}^{(l_1)} \| \leq  \| x_{k} - x_{k}^{(l_2)} \| , \ \text{if} \ l_1 < l_2,   $$ 
where $\| \cdot \|$ is the Euclidean norm in $\mathbb{R}^{d}$. 
Then, for a pre-chosen integer $L$ we use the first $L$ values in $\{\varrho( x_{k}^{(l)} | Y_{1:k}  )\}_{l=1}^N$ to approximate $\varrho(x_{k} | Y_{1:k} )$ as follows 
\begin{equation}\label{Approximation:Prob}
\varrho(x_{k} | Y_{1:k} ) = \sum_{l = 1}^{L} \varrho( x_{k}^{(l)} | Y_{1:k} ) \cdot h_l(x_{k}) , 
\end{equation}
where the weight $h_l(x_{k})$ is defined by
$$h_l(x_{k}) :=  \f{\|x_{k} - x_{k}^{(l)} \| }{\sum_{l = 1}^{L} \|x_{k} - x_{k}^{(l)} \| }. $$
Note that  $\sum_{l = 1}^{L} h_l(x_{k}) = 1$.
From \eqref{Approximation:Prob}, we have
$$
\begin{aligned}
\varrho(x_k | Y_{1:k} ) - p(x_k | Y_{1:k} )  = & \  \sum_{l = 1}^{L} \left( \varrho( x_{k}^{(l)} | Y_{1:k} ) - p( x_{k}^{(l)} | Y_{1:k} ) \right) \cdot h_l(x_{k})  \\ 
& \ + \sum_{l = 1}^{L} p( x_{k}^{(l)} | Y_{1:k} )  \cdot h_l(x_{k})- p( x_{k} | Y_{1:k} ) ,
\end{aligned}
$$
where 
\begin{equation}\label{Approximation:Error}
 \sum_{l = 1}^{L} p( x_{k}^{(l)} | Y_{1:k} )  \cdot h_l(x_{k})- p( x_{k} | Y_{1:k} )  = \sum_{l = 1}^{L} \left( p( x_{k}^{(l)} | Y_{1:k} ) - p( x_{k} | Y_{1:k} ) \right) \cdot h_l(x_{k})
\end{equation}
is the error of the Shepard's interpolation.
We assume that $p(x_{k} | Y_{1:k} )$ has bounded first order derivative. For each pair $p( x_{k} | Y_{1:k} )$ and $p(x_{k}^{(l)} | Y_{1:k} )$ the approximation error $ | p( x_{k}^{(l)} | Y_{1:k} ) - p( x_{k} | Y_{1:k} ) | $ is controlled by the distance $\| x_{k} - x_{k}^{(l)} \|$ and the derivative $p^{\prime}(z | Y_{1:k})$, where $z$ is a point between $x_{k}$ and $x_{k}^{(l)}$.
It is reasonable to assume that  in high probability region of the derivative $p^{\prime}(z | Y_{1:k})$ is large. It's worth pointing out that the random state points generated in this algorithm are concentrated in the high probability region. 
 Thus, if $x_{k}$ lies in the high probability region, the distance $\| x_{k} - x_{k}^{(l)} \|$ is small, which balances the error brought by the large derivative. On the other hand, if $x_{k}$ lies in the low probability region, although the distance $\| x_{k} - x_{k}^{(l)} \|$ is relatively large, the approximation error \eqref{Approximation:Error} is still small due to the small value of the derivative $p^{\prime}(z | Y_{1:k})$.


\subsection{Resampling}\label{Resampling}
Similar to the particle filter method, the above random state points generation suffers from the degeneracy problem for long term simulations, especially for 
high-dimensional problems. After several time steps, the probability density tends to concentrate on a few points which dramatically reduces the number of effective sample points in $\mathcal{H}_k$.

In this work, we propose an occasional resampling procedure to address these problems and rejuvenate the random points cloud. At the time step $k-1$, the resampling procedure takes place after we obtain $\varrho(x_{k-1} | Y_{1:k-1})$, in order to remove the degenerated points in $\mathcal{H}_{k-1}$ using the information provided by $\varrho(x_{k-1} | Y_{1:k-1})$. Specifically, the first step is to develop a degeneracy metric to determine the necessity of doing resampling. To this end, we define the following degenerated subset $\mathcal{S}_{k-1} \subset \mathcal{H}_{k-1}$,
\begin{equation}\label{sss}
\mathcal{S}_{k-1}  = \left\{ x_{k-1}^{i} \big| x_{k-1}^{i} \in \mathcal{H}_{k-1}, \varrho(x_{k-1}^{i} | Y_{1:k-1} ) < \varepsilon \right\},
\end{equation}
where $\varepsilon > 0$ is a user-defined threshold. We also define
\[
\mathcal{J}(\mathcal{S}_{k-1}):= \{i = 1, \ldots, N | x_{k-1}^i \in \mathcal{S}_{k-1}\}
\]
to be the index set of $\mathcal{S}_{k-1}$.
Then, the degeneracy of $\mathcal{H}_{k-1}$ can be measured by the ratio $\#(\mathcal{S}_{k-1}) / \#(\mathcal{H}_{k-1}) \in [0,1]$, where $\#(\cdot)$ denotes the number of points in a set. If the ratio is smaller than a threshold $\tau \in [0, 1]$, then we will skip the resampling step and propagate $\mathcal{H}_{k-1}$ to get $\mathcal{H}_{k}$; otherwise, the set $\mathcal{H}_{k-1}$ is considered degenerated, and the resampling procedure is needed. 


In resampling, instead of propagating $\mathcal{H}_{k-1} $ to $\mathcal{H}_{k} $, we aim at constructing an intermediate point set, denoted by $\mathcal{H}_{k-\f{1}{2}} := \{ x^i_{k-\f{1}{2}}\}_{i=1}^N$ and propagate $\mathcal{H}_{k-\f{1}{2}}$ through the state model \eqref{eq:State} to obtain $\mathcal{H}_{k}$. According to the definition of $\mathcal{S}_{k-1}$ in \eqref{sss}, we consider the state points in $\mathcal{H}_{k-1} \backslash \mathcal{S}_{k-1}$ are in the statistically significant region of $\varrho(x_{k-1} | Y_{1:k-1} ) $, so that we first put those points in $\mathcal{H}_{k-\f{1}{2}}$, i.e.,
\[
x_{k-\frac{1}{2}}^i = x_{k-1}^i\; \text{ for } \; i \notin \mathcal{J}(\mathcal{S}_{k-1}).
\]
For the state points in $\mathcal{S}_{k-1}$, we replace them by generating new samples from $\varrho(x_{k-1} | Y_{1:k-1} )$ using the importance sampling \cite{Budhiraja-Survey}, i.e.,
\[
x_{k-\frac{1}{2}}^i \sim \varrho(x_{k-1} | Y_{1:k-1} )\; \text{ for } \; i \in \mathcal{J}(\mathcal{S}_{k-1}).
\]
%
%
%
As a result, the resampling procedure helps us remove the state points with low probabilities, and makes the state point set $\mathcal{H}_{k}$ concentrated in the high probability region of the posterior PDF $\varrho(x_{k-1} | Y_{1:k-1} )$ at each time step. 
%


\subsection{Summary of the algorithm}
Finally, we summarize the entire meshfree implicit filter algorithm introduced in \S \ref{Prediction}-\S \ref{Resampling} in 
Algorithm 1 below.
%

\begin{table}[h!]
\begin{tabular}{p{0.95\textwidth}}
\hline\noalign{\smallskip}
{\bf Algorithm 1}: {\em The meshfree implicit filter algorithm}\\
\noalign
{\smallskip}\hline
\noalign{\smallskip}
\vspace{-0.3cm}
\begin{spacing}{1.1}
\begin{algorithmic}\label{algorithm2}
\vspace{0.0cm}
\State {\bf Initialization}: set the number of samples $M$ for estimating $\mathbb{E}_w[\cdot]$,
 the number of state points $N$, the resampling thresholds $\varepsilon$ and $\tau$
 \vspace{0.1cm}
\While{$k = 1, 2, \cdots,$}
\State Compute the ratio ${\#(\mathcal{S}_{k-1})}/{\#(\mathcal{H}_{k-1})}$
 \vspace{0.1cm}
 \If{${\#(\mathcal{S}_{k-1})}/{\#(\mathcal{H}_{k-1})} < \tau$}
 \vspace{0.1cm}
 \State Propagate $\mathcal{H}_{k-1}$ through the state model \eqref{eq:State} to obtain $\mathcal{H}_{k}$
 \vspace{0.1cm}
 \Else 
 \State Resample and construct the intermediate state set $\mathcal{H}_{k-\f{1}{2}}$
  \vspace{0.1cm}
 \State Propagate $\mathcal{H}_{k-\f{1}{2}}$ through the state model \eqref{eq:State} to obtain $\mathcal{H}_{k}$
   \vspace{0.1cm}
 \EndIf
\vspace{0.1cm}
\State {\bf Prediction}: solve $\varrho(x_{k}| Y_{1:k-1} )$ using \eqref{scheme:a}, at each point in $\mathcal{H}_k$
\vspace{0.1cm}
\State {\bf Update}: solve $\varrho(x_{k}| Y_{1:k} )$ using \eqref{scheme:b} and \eqref{Approximation:Prob}
\vspace{0.1cm}
\EndWhile
\end{algorithmic}
\vspace{-0.9cm}
\end{spacing}
\\
\hline
\end{tabular}
\end{table}

\vspace{0.1in}

\section{Numerical experiments}\label{sec:ex}
In this section, we present two numerical examples to examine the performance of our meshfree implicit filter method. In Example 1, we use a two dimensional nonlinear filtering problem to show the distributions of the random points $\mathcal{H}_k$. In Example 2, we solve a three dimensional bearing-only tracking problem, which is a six dimensional nonlinear filtering problem. For this higher dimensional problem, we compare the accuracy and efficiency of our meshfree implicit filter method with the extended Kalman filter and the particle filter. 
\subsection*{Example 1}\label{Ex1} 
In this example, we consider the two dimensional noise perturbed tumoral growth model \cite{2D_PopulationModel}

\begin{equation}\label{Eq:Ex1:Cont}
d \bm{X}_t = F(\bm{X}_t)  dt + \bm{\sigma} \cdot dW_t,
\end{equation}
where $W_t$ is a two dimensional standard Brownian motion and $\bm{\sigma} = (0.01, 0.01)^T$.  The state process $\bm{X}_t = (X_t^1, X_t^2)^T$ is  a two dimensional vector, $F(\bm{X}_t) := ( f_1(\bm{X}_t), f_2(\bm{X}_t) )^T $ is defined as
$$ f_1(\bm{X}_t) = \alpha_1 X_t^1 \cdot \ln(\f{X_t^2}{X_t^1}) $$
and
$$ f_2(\bm{X}_t) = \alpha_2 X_t^1 - \alpha_3 X_t^2 \cdot (X_t^1)^{\f{2}{3}}. $$
Here, $f_1$ models the Gompertzian growth rate of the tumor and $f_2$ gives the degree of vascularization of the tumor which is also called `` angiogenic capacity''.

To approximate the state variables, we discretize the dynamic system \eqref{Eq:Ex1:Cont} in time and obtain a discrete state model
\begin{equation}\label{Eq:Ex1:Discrete}
\bm{X}_k = F(\bm{X}_{k-1}) \cdot \Delta + \bm{\sigma} \cdot \bm{\omega}_{k-1}.
\end{equation}
Here, $\bm{\omega}_k$ is a two dimensional zero mean Gaussian white noise process with covariance $Q = I\Delta$, where $I$ is the $2 \times 2$ identity matrix and $\Delta$ is the time partition stepsize.

The measurement of the state model is given by
$$  Y_k =\left( X_k^1, X_k^2 \right)^T + \bm{R} \cdot \bm{v}_k,   $$
where $\bm{v}_k$ is a two dimensional zero mean Gaussian white noise process with covariance $\Lambda = I \Delta$, $I$ is a $2 \times 2$ identity matrix and $\bm{R} = (0.1, 0.1)^T$.

In the numerical experiment, we use uniform time partition with stepsize $\Delta = 0.2$ and simulate the state process for $K = 40$ with initial state $\bm{X}_0 = ( 0.8, 0.3 )^T$ and parameters $\alpha_1 = 1$, $\alpha_2 = 0.2$, $\alpha_3 = 0.2$. At time step $k = 0$, we initialize the prior PDF $p_0$ by $N(\tilde{\bm{X}}_0, \Sigma)$, where $\tilde{\bm{X}}_0 = (0.78, 0.32)^T$ and
\begin{align} \label{Ex1:Ini_sigma}
\Sigma = \left(
\begin{array}{cc}
0.05^2 & 0 \\
0& 0.1^2 \\
\end{array}\right).
\end{align}
\begin{figure}[ht!]
\begin{center}
  \includegraphics[scale = 0.5]{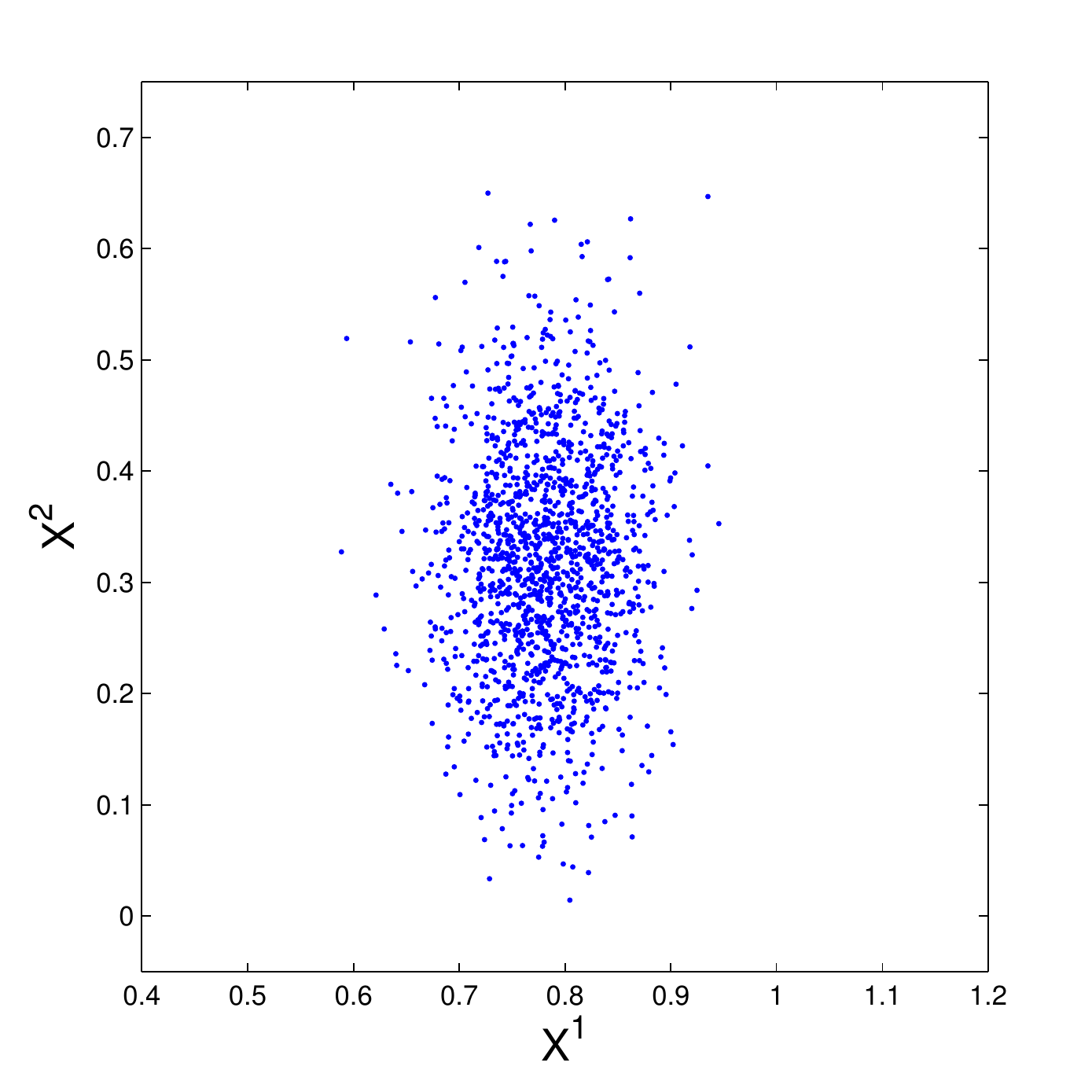}
\end{center}
\caption{Example 1 : Initial random state space points $\mathcal{H}_0$}\label{scatter_0}
\end{figure}
\begin{figure}[ht!]
\begin{center}
\subfloat[ k = 1 ]{\includegraphics[scale = 0.4]{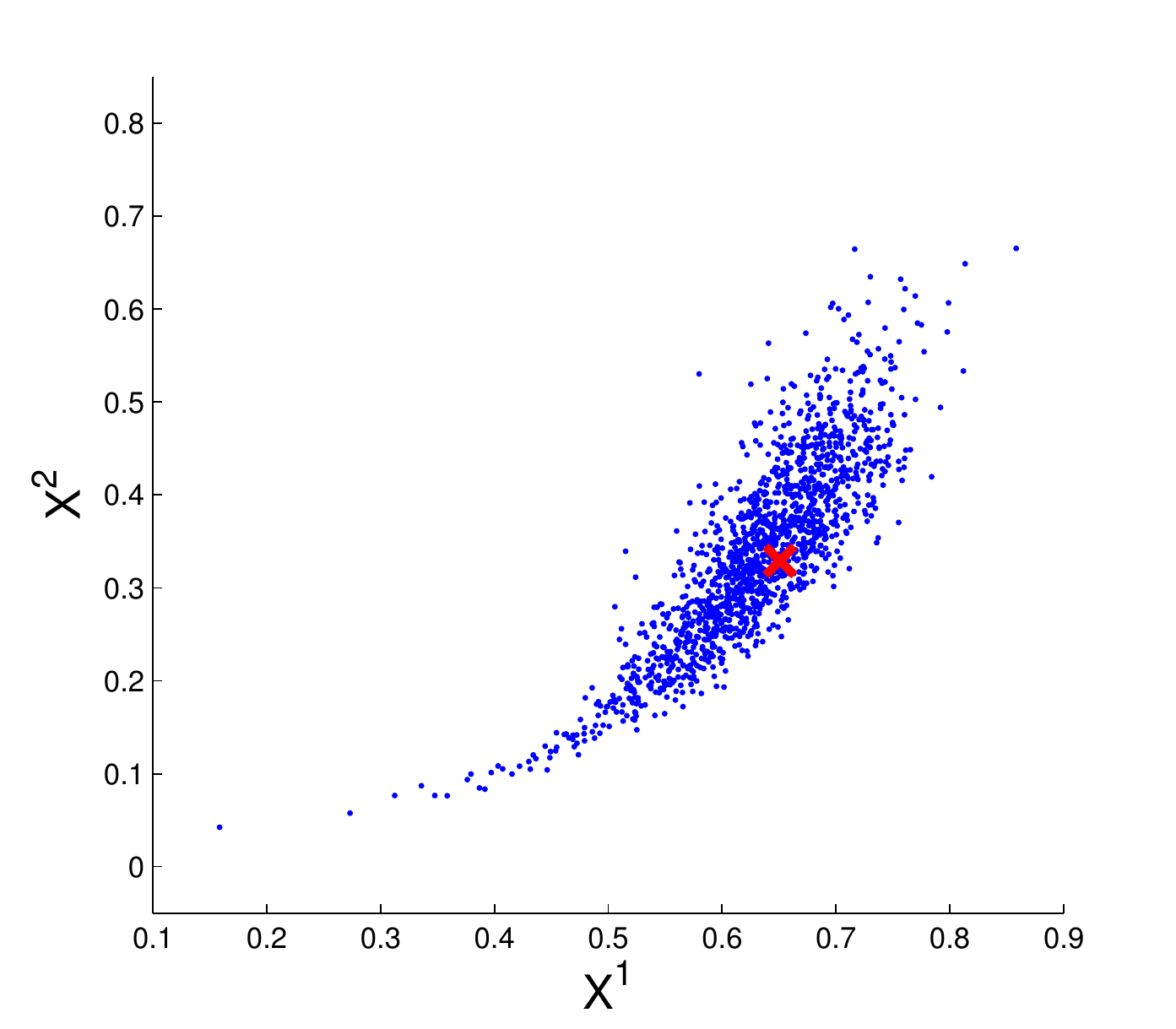}}\label{scatter_1}
\subfloat[ k = 2 ]{\includegraphics[scale = 0.4]{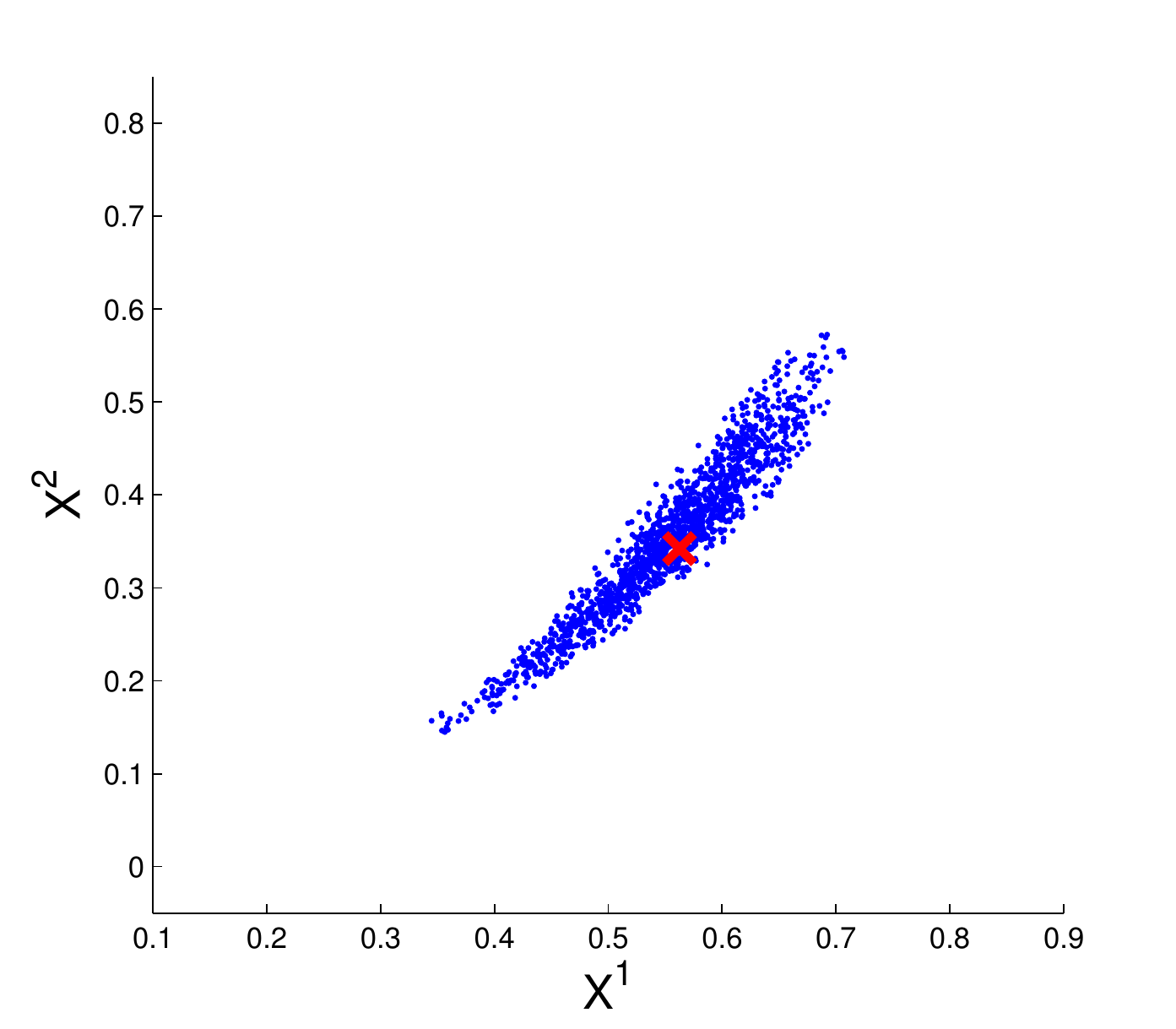}}\label{scatter_2}\\
\subfloat[ k = 3 ]{\includegraphics[scale = 0.4]{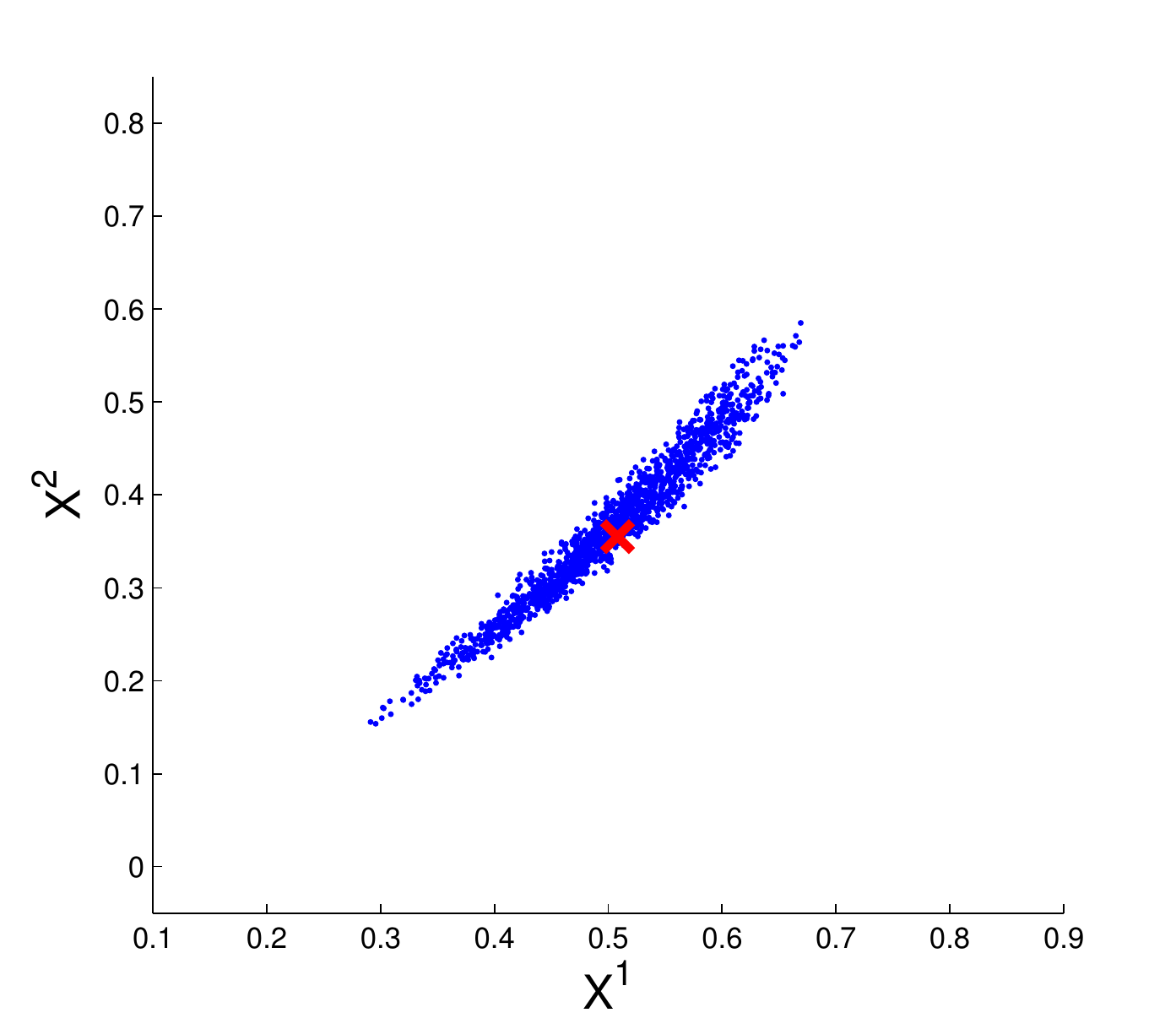}}\label{scatter_3}
\subfloat[ k = 10 ]{\includegraphics[scale = 0.4]{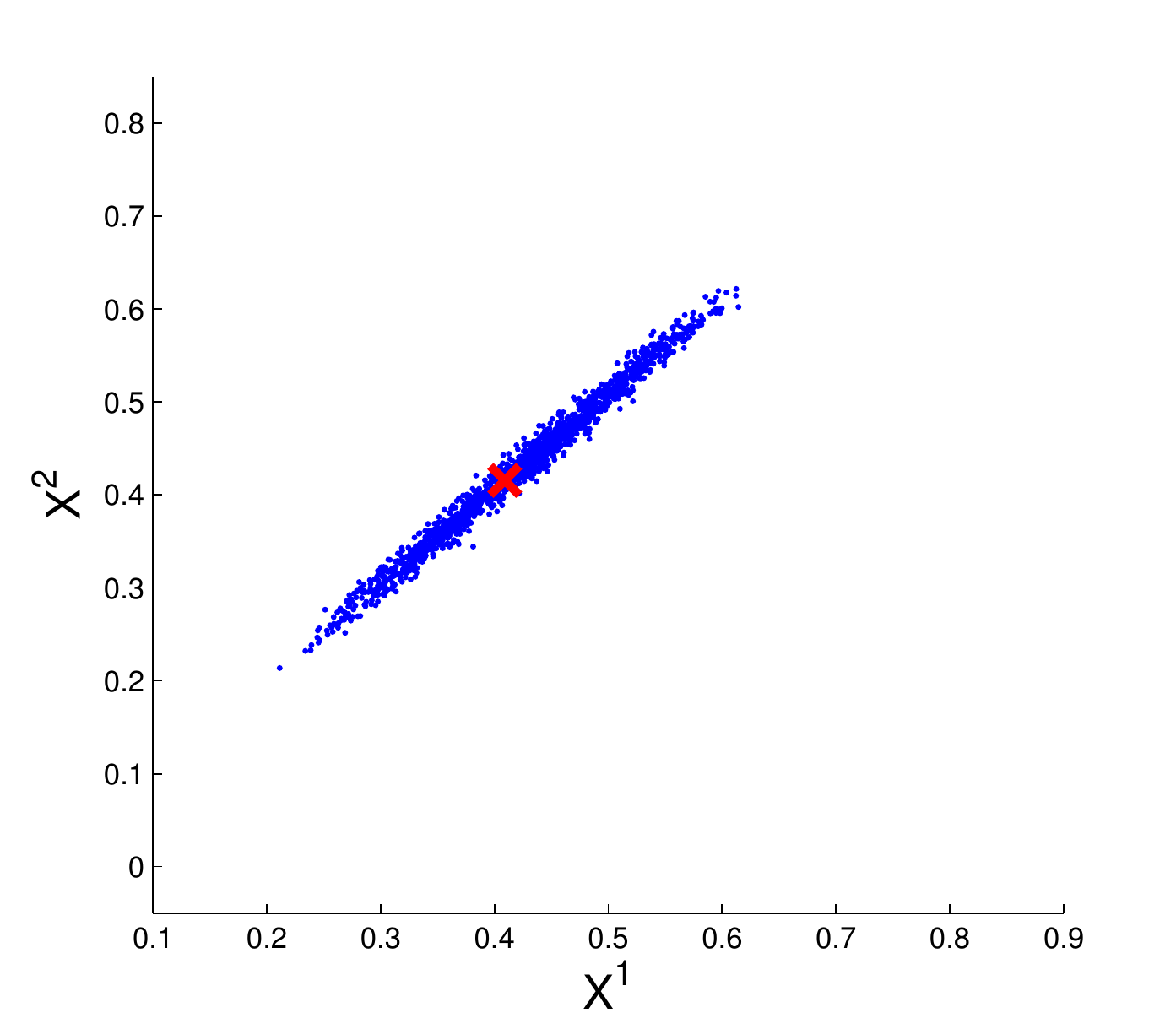}}\label{scatter_10}\\
\subfloat[ k = 20 ]{\includegraphics[scale = 0.4]{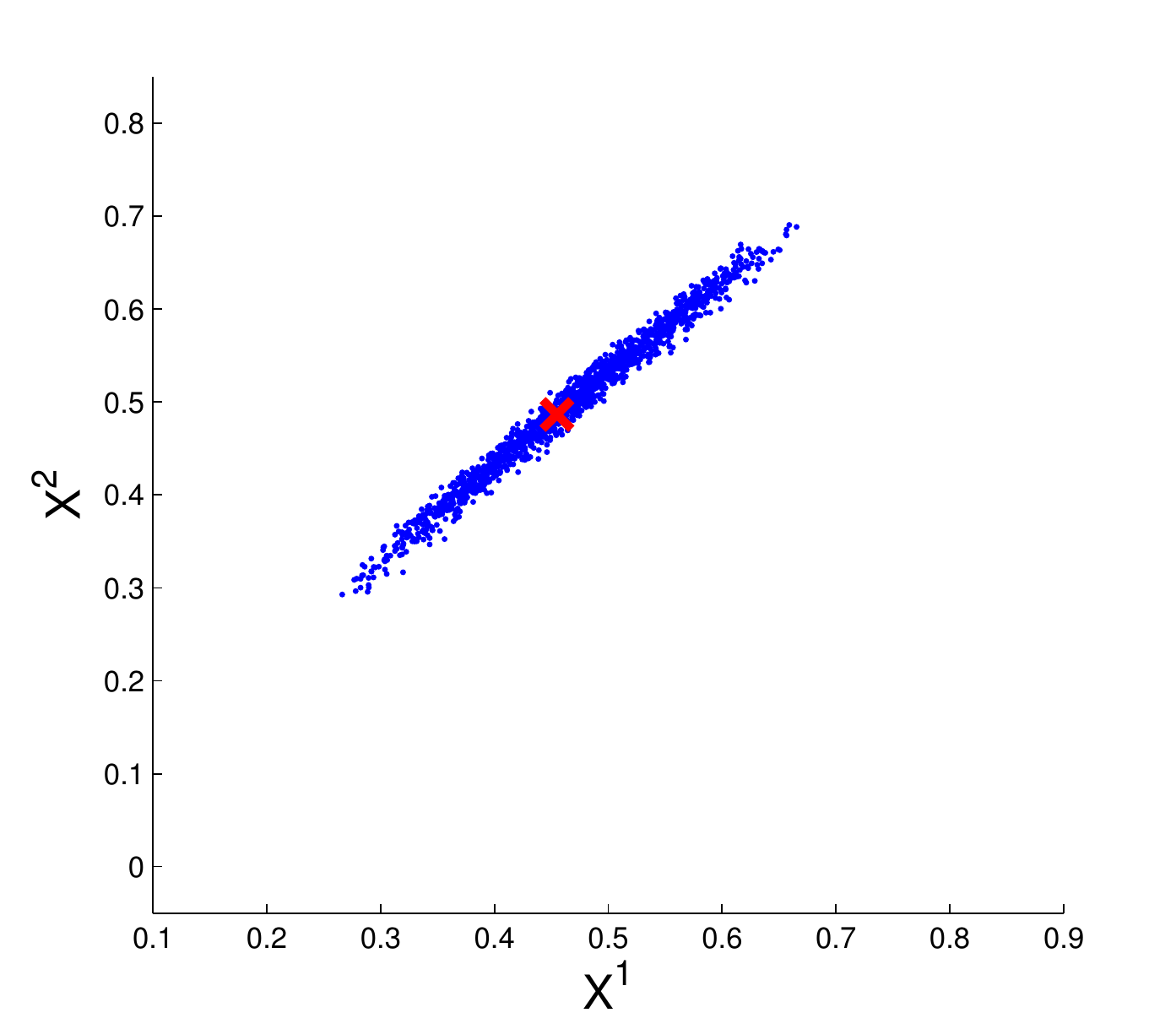}}\label{scatter_20}
\subfloat[ k = 40 ]{\includegraphics[scale = 0.4]{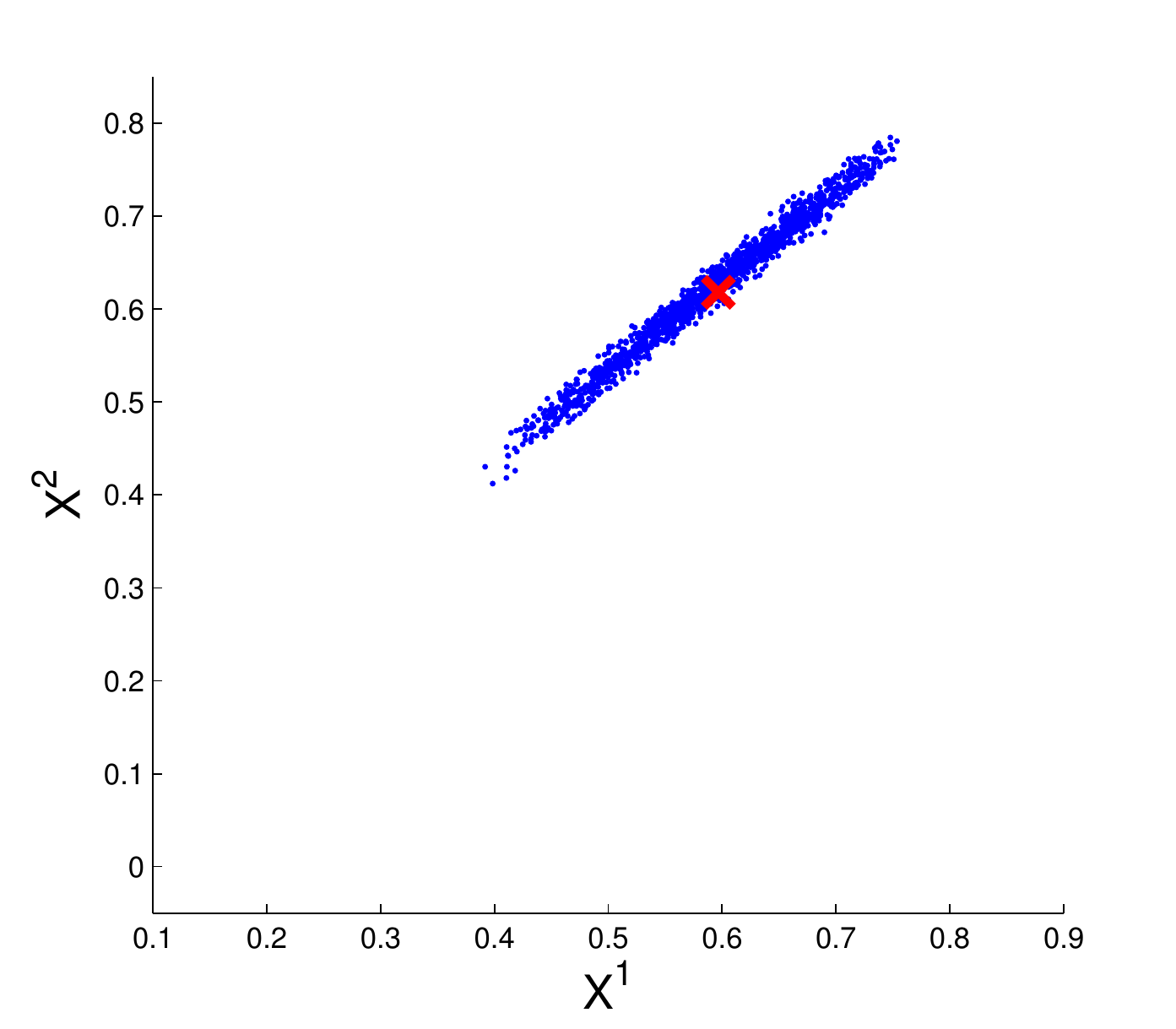}}\label{scatter_40}\\
\end{center}
\caption{Example 1: State space points $\mathcal{H}_k$ at time steps $k = 1, 2, 3, 10, 20, 40$. }\label{2D_Grids}
\end{figure}

In Figure \ref{scatter_0}, we plot $1500$ random samples generated from the initial PDF $p_0$, which are our initial random points $\mathcal{H}_0$.  Figure \ref{2D_Grids} illustrates the behavior of random state points $\mathcal{H}_k$ at time steps $k = 1, 2, 3, 10, 20, 40$, respectively. In Figure \ref{2D_Grids}, the blue dots in each figure plot the random state points obtained by using the dynamic state points generation method introduced in Section \ref{Algorithm} and the red cross in each figure gives the true state $\bm{X}_k$ at the corresponding time step. From the figures we can see that all the points are moving according to the state model and are concentrated around the true state.

\begin{figure}[ht!]
\begin{center}
\subfloat[ Simulation: $X^1$ ]{\includegraphics[scale = 0.42]{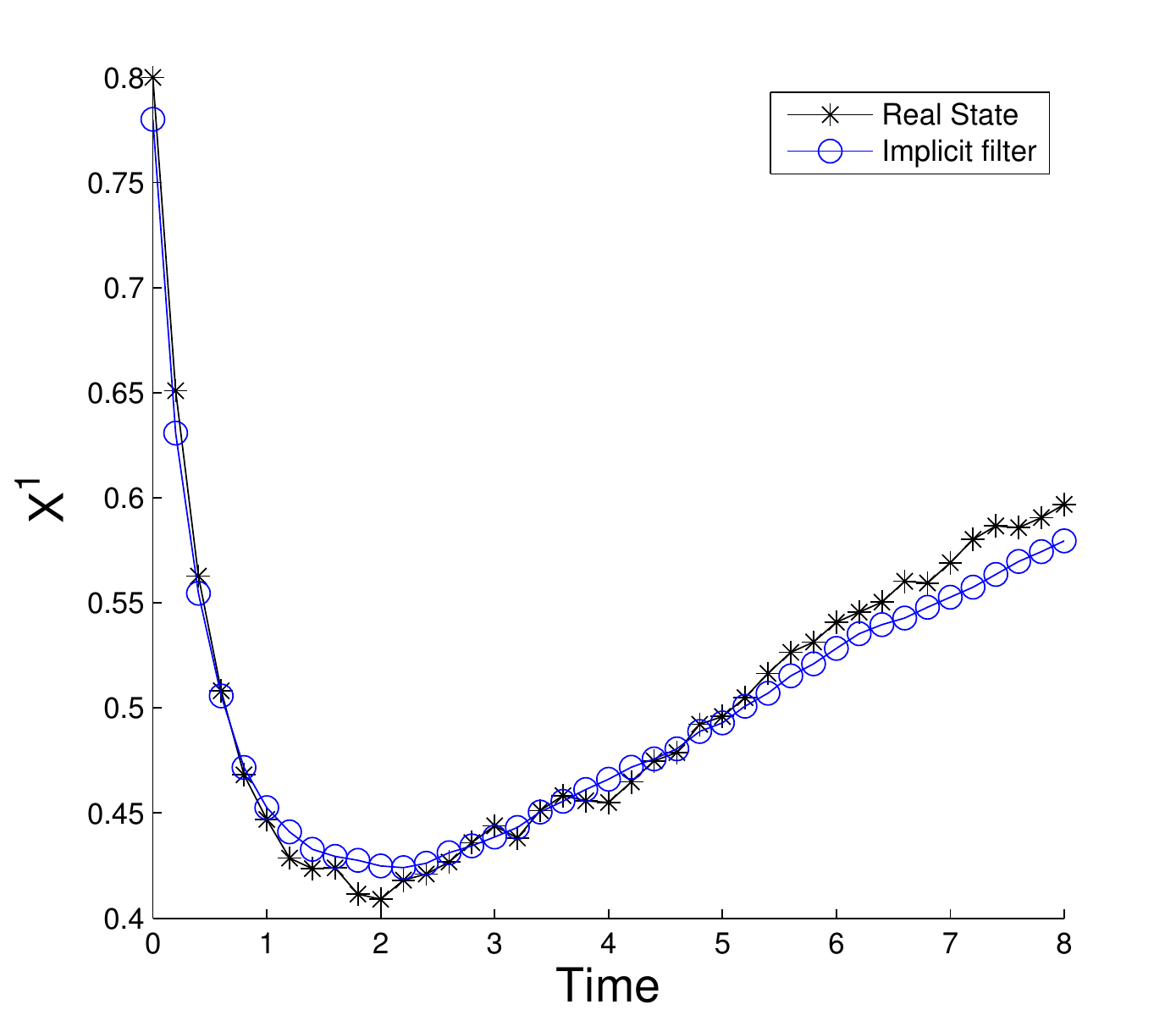}}\label{2D_Simulation_X1}
\subfloat[ Simulation: $X^2$ ]{\includegraphics[scale = 0.42]{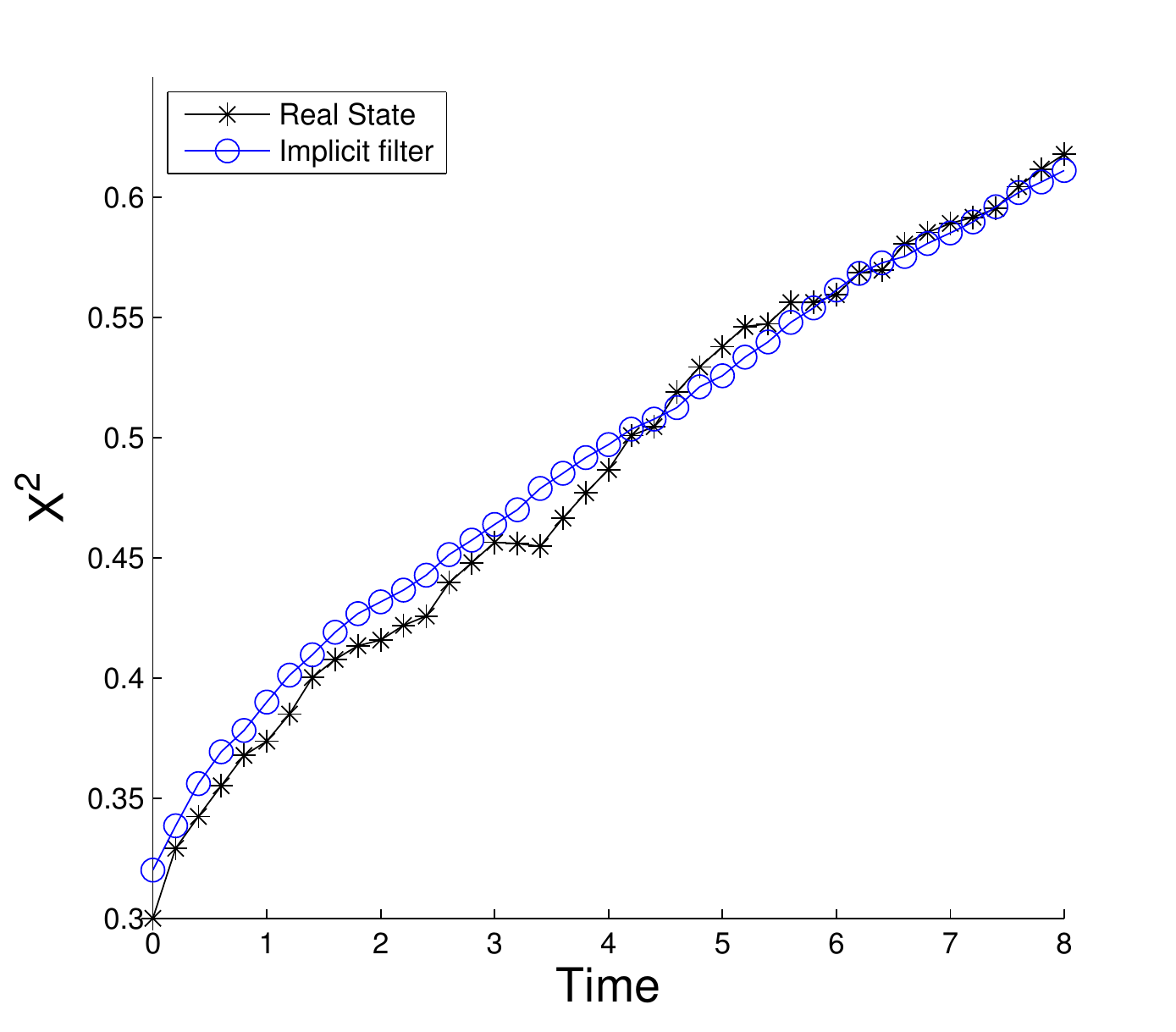}}\label{2D_Simulation_X2}
\end{center}
\caption{Example 1: Simulation of the tumoral growth states }\label{2D_Simulation}
\end{figure}

To present the accuracy of the algorithm, we show the simulation of the tumoral growth states in Figure \ref{2D_Simulation}. The black curves are the true  $X^1$ and $X^2$ coordinate values of the tumoral growth states, respectively. The blue curves show the simulated states obtained by using the meshfree implicit filter method.

\subsection*{Example 2}
\begin{figure}[ht!]
\begin{center}
  \includegraphics[scale = 0.48]{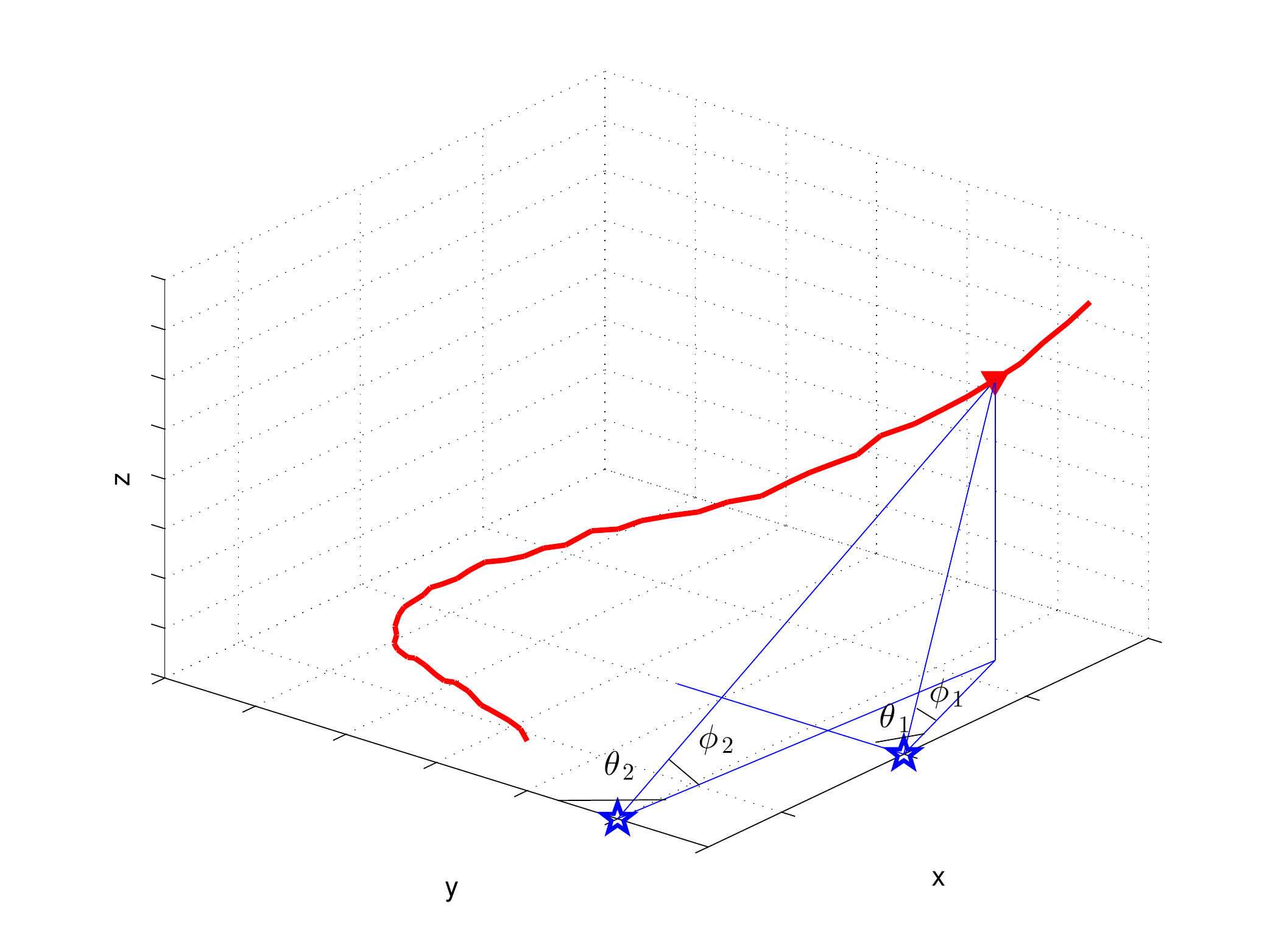}
\end{center}
\caption{Example 2: Bearing-only Tracking in 3-D}\label{Model_6D}
\end{figure}

In this example, we study a six dimensional target tracking problem. 
In Figure \ref{Model_6D},  the  target,  denoted by the red  line,  moves in the three dimensional space and  two platforms on the ground,  denoted by pentagons,   take angular  observations of the moving target.

\begin{figure}[ht!]
\begin{center}
\includegraphics[scale = 0.5]{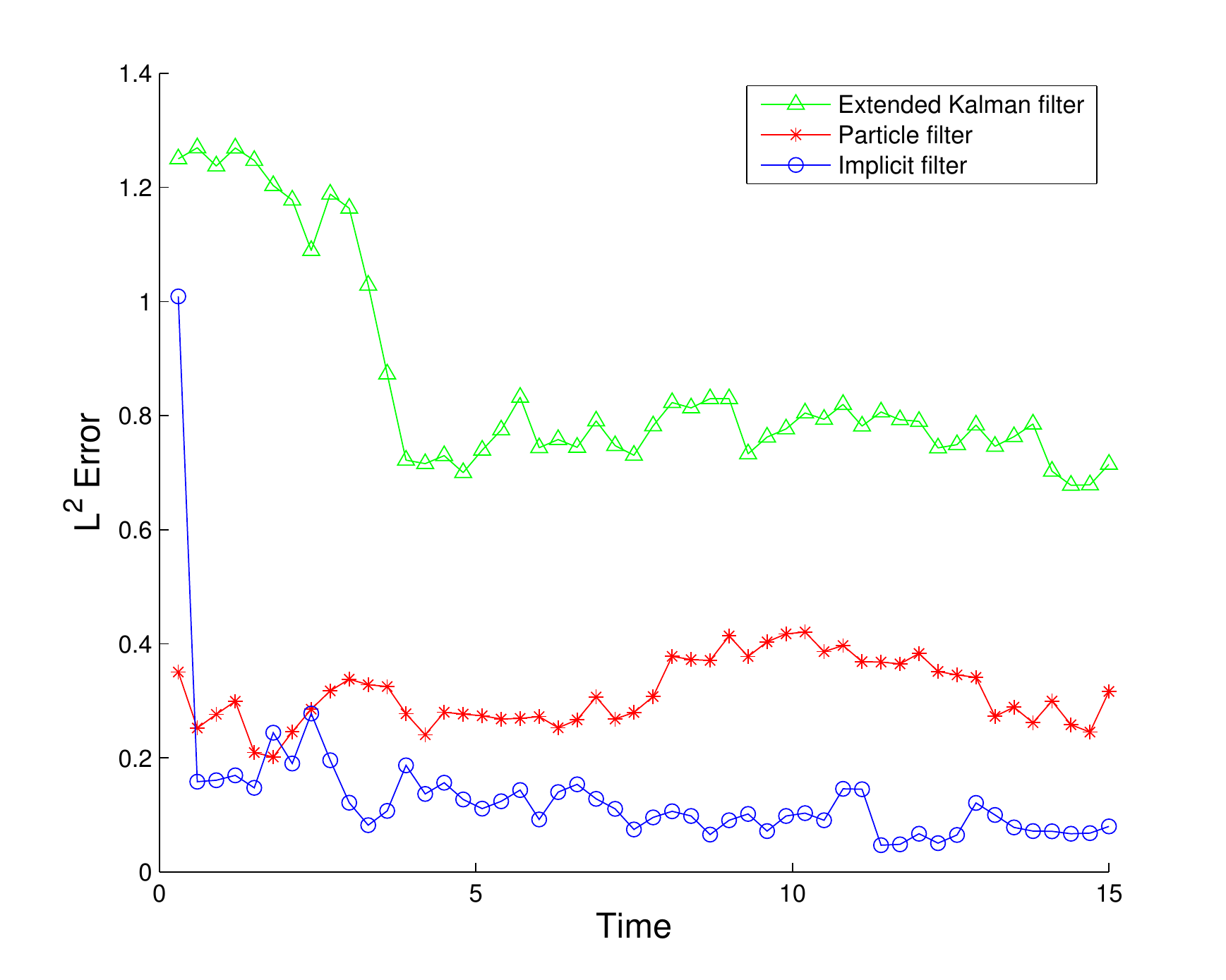}
\end{center}
\caption{Example 2 : Comparison of $L^2$ error. }\label{6D_L2}
\end{figure}

The state process $\bm{X}_k = (X^1_k, X^2_k, X^3_k, X^4_k, X^5_k, X^6_k)^T$ is described by the following dynamic model
\begin{equation}\label{Eq:Ex3:State}
\bm{X}_k = f(\bm{X}_{k-1}) + \bm{\sigma} \cdot \bm{\omega}_{k-1},
\end{equation}
where  $(X^1, X^2, X^3)$ describes the position of the moving target which is controlled by parameters $(X^4, X^5, X^6)$.  The system noise $\bm{\omega_{k}} = (\omega^1_{k}, \omega^2_{k}, \omega^3_{k}, \omega^4_{k}, \omega^5_{k}, \omega^6_{k})^T$ is a zero mean Gaussian white noise process with covariance $Q \doteq I \Delta $, $I$ is the $6 \times 6$ identity matrix and $\Delta$ is a given time period, $ \bm{\sigma} = (0.1, 0.1, 0.1, 0.01, 0.01, 0.01 )^T$ is a constant vector and $f$ is given by
$$
\renewcommand{\arraystretch}{1.5}
f(\bm{X}_k) = \left(
\begin{array}{ccc}
X_{k-1}^1 + X_{k-1}^4 \Delta  \\
X_{k-1}^2 + \sin( \alpha X_{k-1}^5) \Delta  \\
X_{k-1}^3 + (X_{k-1}^6)^2 \Delta \\
X_{k-1}^4 + v_1 \Delta \\
X_{k-1}^5 + v_2 \Delta \\
X_{k-1}^6 + v_3 \Delta
\end{array}\right).
$$
The measurements $\bm{Y}_k$ of the state process from the two locations are given by 
$$
\renewcommand{\arraystretch}{1.5}
\bm{Y}_k = \left(
\begin{array}{ccc}
 \arctan\left( \f{X_k^3}{ \sqrt{( X_k^1 - a_1)^2 + ( X_k^2 - b_1 )^2} } \right) \\
 \arctan\left( \f{X_k^3}{ \sqrt{( X_k^1 - a_2)^2 + ( X_k^2 - b_2 )^2} } \right) \\
 \arctan\left( \f{X_k^1 - a_1}{ X_k^2 - b_1 } \right) \\
 \arctan\left( \f{X_k^1 - a_2}{ X_k^2 - b_2 } \right)
\end{array}\right) + \bm{R} \bm{v}_k,
$$
where $\bm{v}_k$ is a 4 dimensional zero mean Gaussian white noise process with covariance $\Lambda = I \Delta$, $I$ is a $4 \times 4$ identity matrix, $\bm{R} = (0.6, 0.6, 0.6, 0.6)^T$, $(a_1, b_1)$ and $(a_2, b_2)$ are locations of two observers.
\begin{figure}[ht!]
\begin{center}
\subfloat[]{\includegraphics[scale = 0.48]{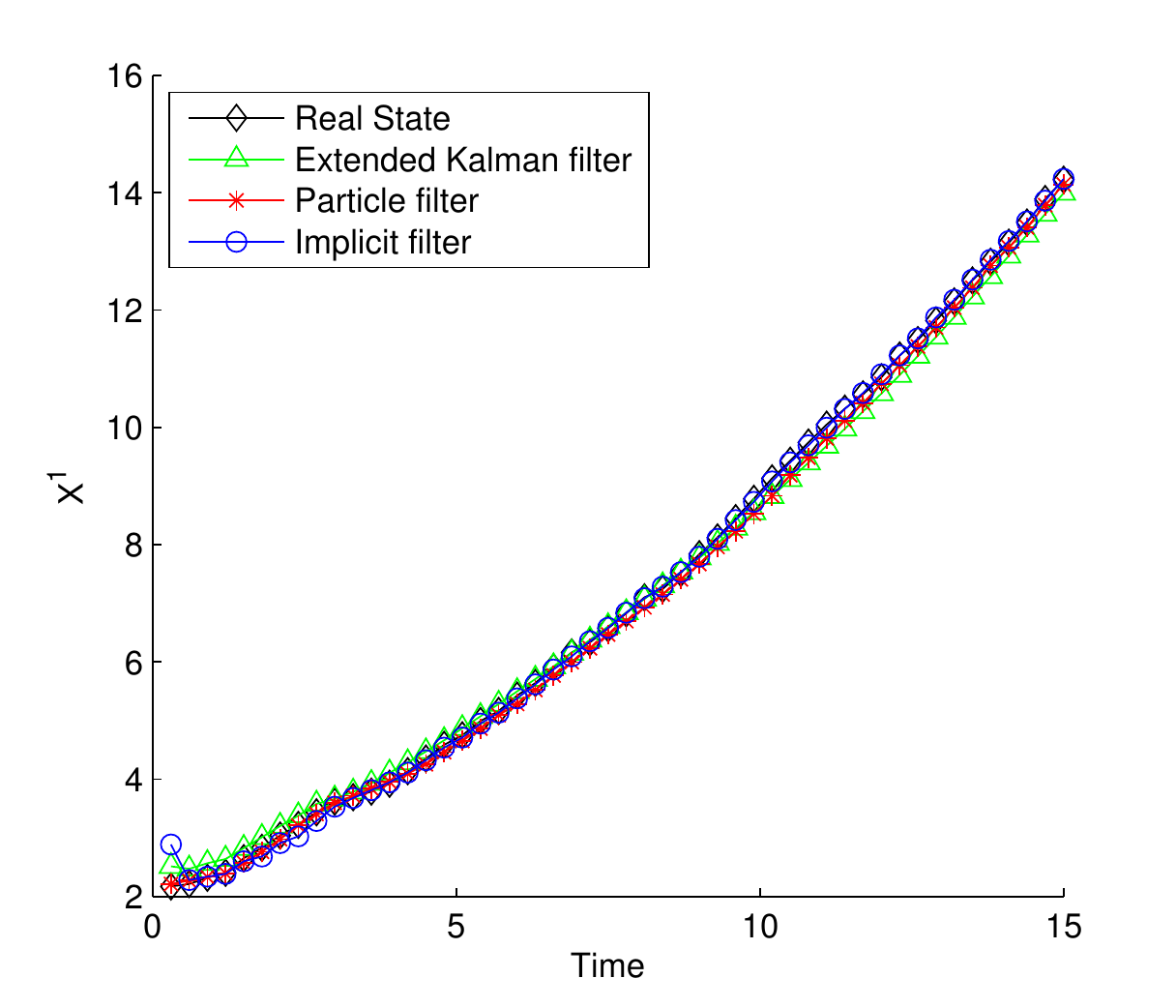}}\label{Com_6D_X1}
\subfloat[]{\includegraphics[scale = 0.48]{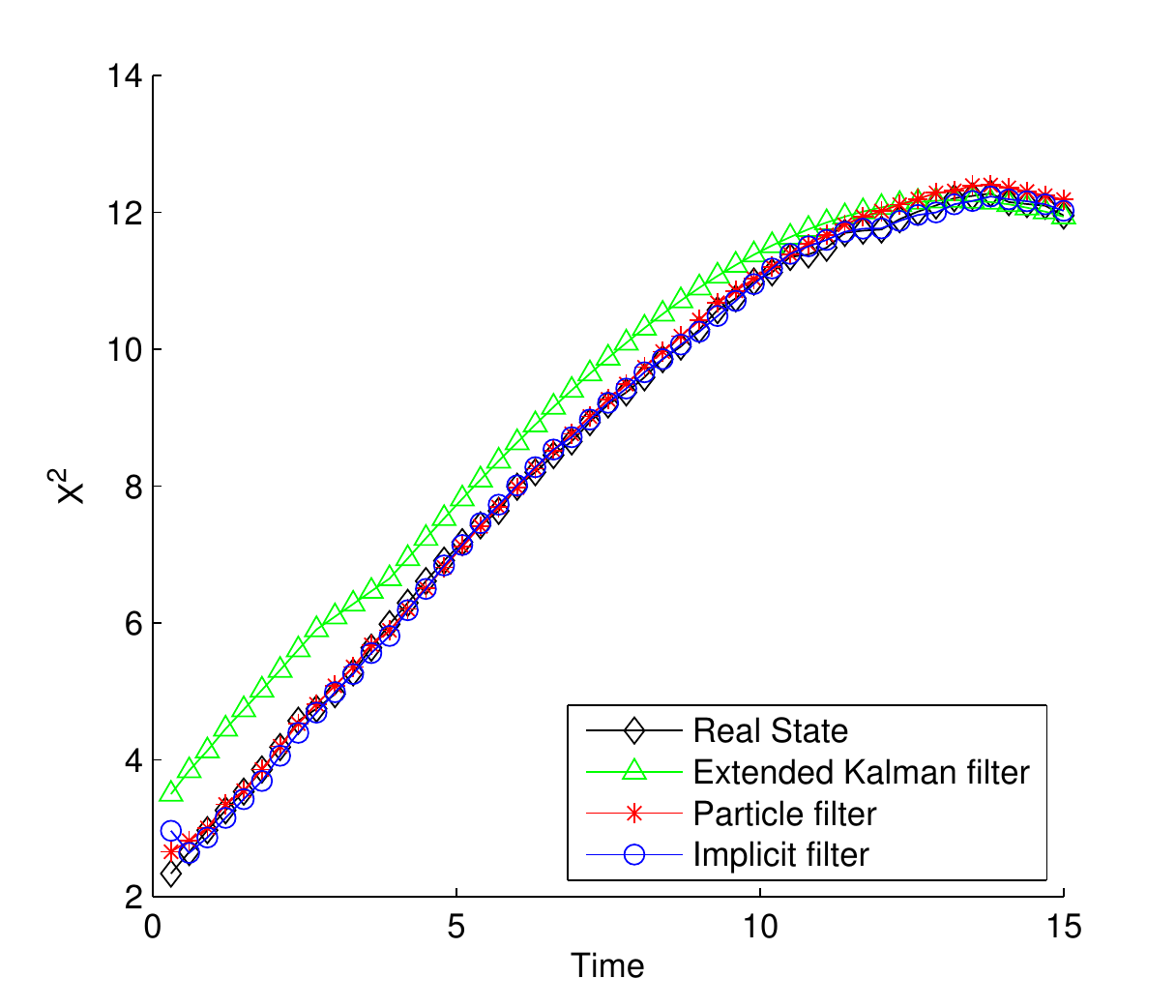}}\label{Com_6D_X2}\\
\subfloat[]{\includegraphics[scale = 0.48]{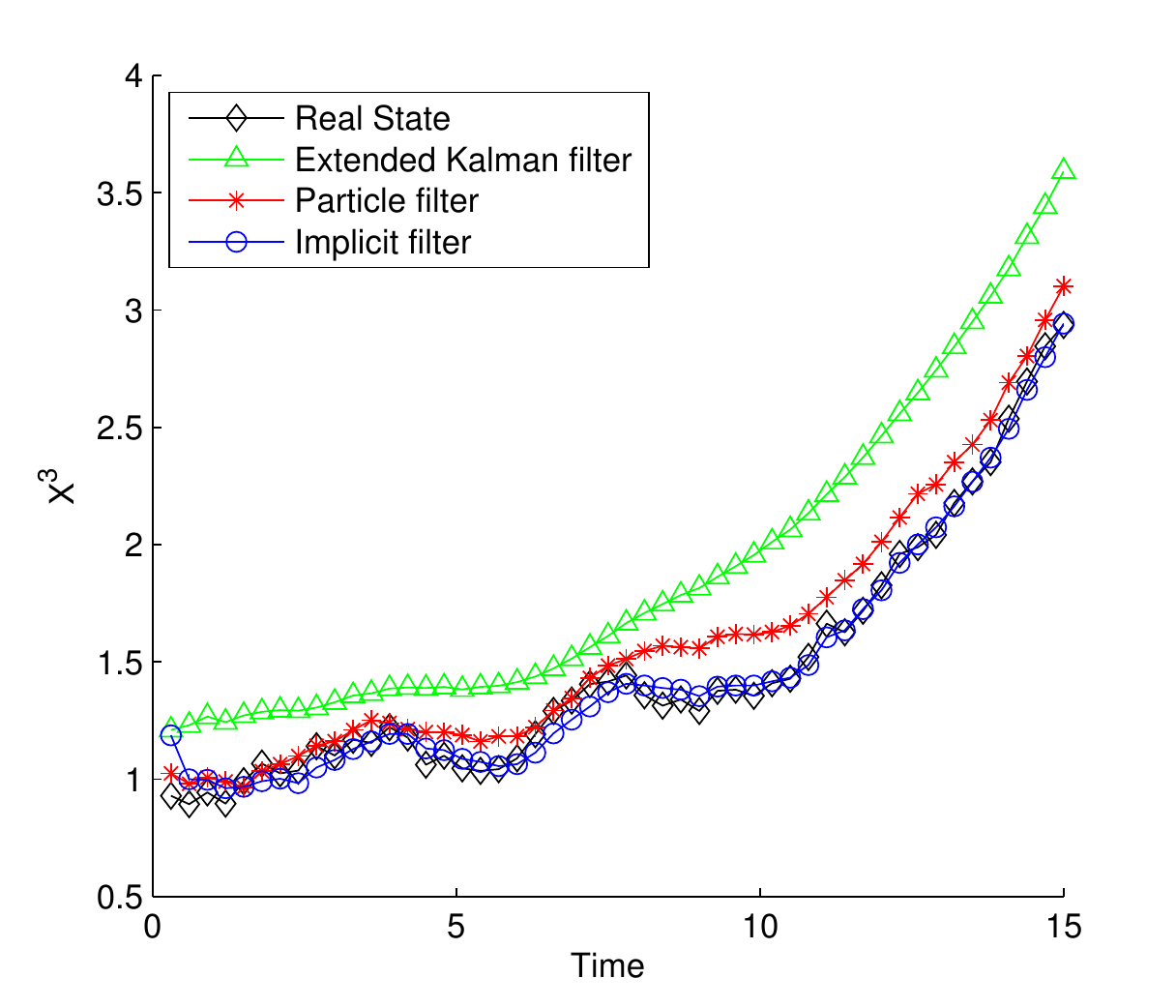}}\label{Com_6D_X3}
\subfloat[]{\includegraphics[scale = 0.48]{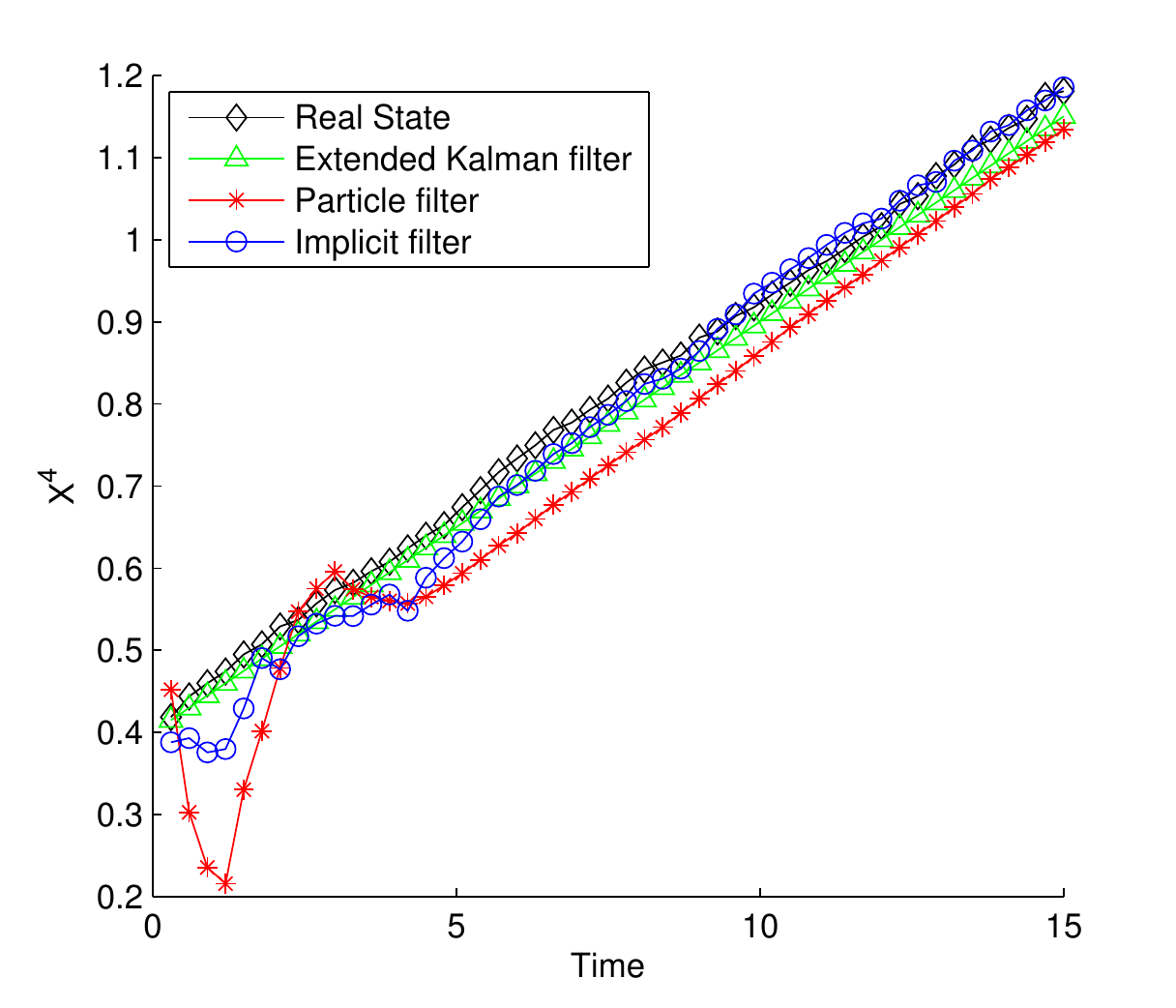}}\label{Com_6D_X4}\\
\subfloat[]{\includegraphics[scale = 0.48]{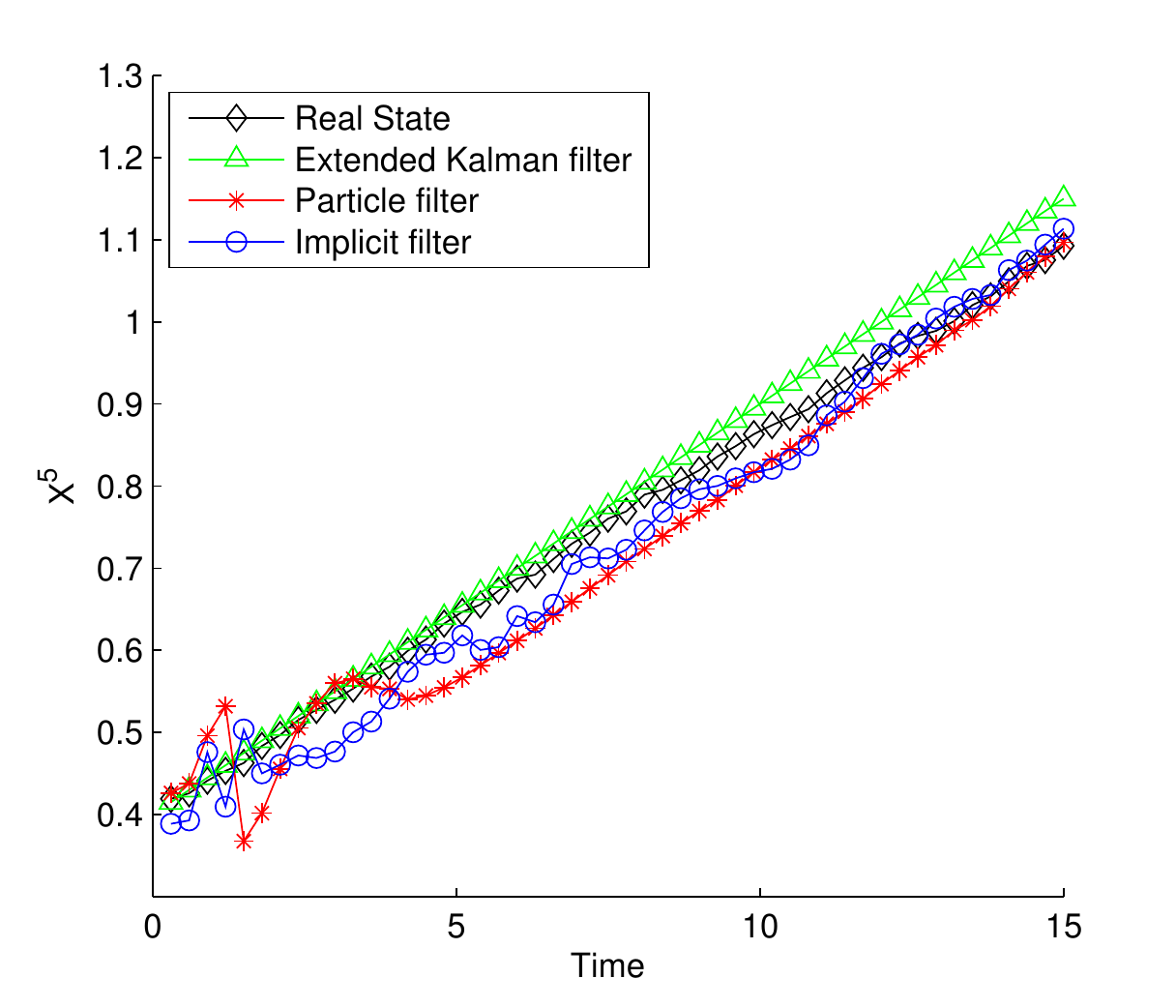}}\label{Com_6D_X5}
\subfloat[]{\includegraphics[scale = 0.48]{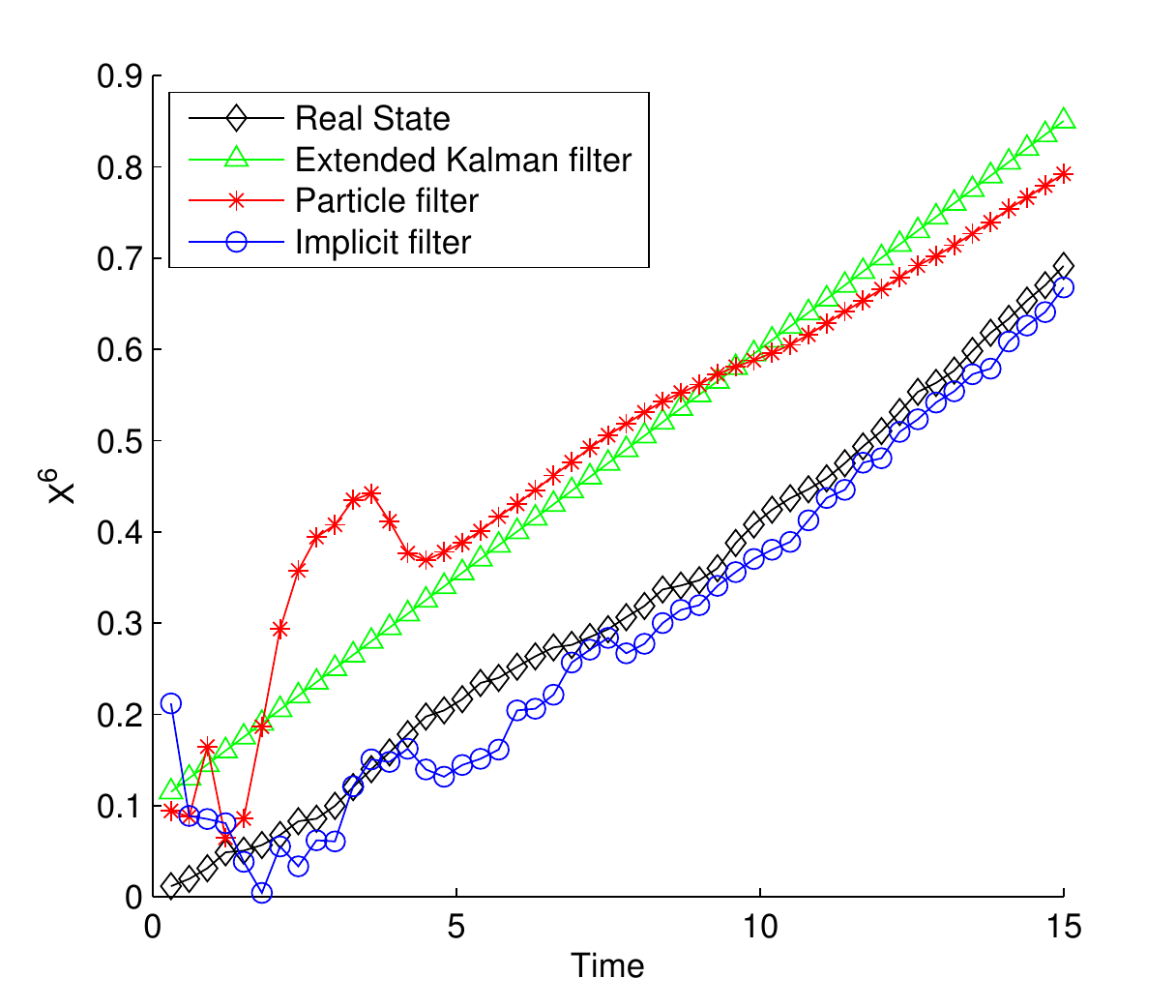}}\label{Com_6D_X6}\\
\end{center}
\caption{Example 2 : Comparison of estimated states. (a) Shows the comparison on $X^1$ direction. (b) Shows the comparison on $X^2$ direction. (c) Shows the comparison on $X^3$ direction. (d) Shows the comparison on $X^4$ direction. (e) Shows the comparison on $X^5$ direction. (f) Shows the comparison on $X^6$ direction. }\label{6D_State_Comparison}
\end{figure}
We choose $\Delta = 0.3$, $\alpha = 3$, $v_1 = v_2 = v_3 = 0.05$ . Also, we assume that platforms are located at $(a_1, b_1) = (16, 6)$, $(a_2, b_2) = (8, 15)$ and the initial sate is given by a Gaussian $N( \bm{X}_0, \Sigma)$ where  $\bm{X}_0 = ( 2, 2, 1, 0.4, 0.4, 0 )^T$ and

\begin{align*}
\Sigma = \left(
  \begin{array}{cccccc}
    1^2 & 0 & 0 & 0 & 0 & 0 \\
    0 & 1^2 & 0 & 0 & 0 & 0 \\
    0 & 0 & 1^2 & 0 & 0 & 0 \\
    0 & 0 & 0 & 0.2^2 & 0 & 0 \\
    0 & 0 & 0 & 0 & 0.2^2 & 0 \\
    0 & 0 & 0 & 0 & 0 & 0.2^2 
  \end{array}
\right).
\end{align*}
The target will be observed over the time period $0\leq t \leq 15$.  In the numerical experiments, we  compare the performance of our meshfree implicit filter with the extended Kalman filter and the particle filter. In particular, we compare the estimated mean values of the states process along each dimension in Figure \ref{6D_State_Comparison}.
In the particle filter method, we choose $15,000$ particles. In the meshfree implicit filter method, we choose the number of state points to be $N = 4,000$ and the number of random samples in the implicit filter Monte Carlo simulation to be $M = 6$. The black curves in Figure \ref{6D_State_Comparison} show the real states process along each direction, the green curves give the estimated means obtained by the extended Kalman filter method, the red curves give the estimated means obtained by the particle filter method,  and the blue curves give the estimated means obtained by the meshfree implicit filter. We also plot the $L^2$ error $err_k$ corresponding to all three methods in figure \ref{6D_L2}.
As we can see from figure \ref{6D_State_Comparison} and \ref{6D_L2}, the implicit filter and the particle filter are much more accurate than the extended Kalman filter and the implicit filter is the most accurate approximation in this experiment.

To further compare the efficiency between the meshfree implicit filter and the particle filter, we repeat the above experiment over $50$ realizations and show the average CPU time and the corresponding global root mean square error $err_G$ defined by
$$err_G^2 = \f{1}{50}\f{1}{K} \sum_{j=1}^{50}\sum_{k = 1}^{K} err_k(j)^2 $$
where $err_k(j)$ is the $L^2$ error of the $j$-realization at time step $k$.
In table \ref{efficiency},  we can see that with $15,000$ particles, the CPU time of the particle filter method is comparable to that of the implicit filter with $4,000$ random state points, but  the global RMSE of the particle filter is more than doubled the RMSE of the implicit filter. With $25,000$ particles, the particle filter method achieves an accuracy comparable to the implicit filter, but at a significantly higher cost. 

\renewcommand{\arraystretch}{1.25}
\begin{table} 
\leftmargin=6pc \caption{Example 2: Efficiency comparison} \label{efficiency} \small
\begin{center}
\begin{tabular}{|c|c|c|}
 \hline  Methods &  CPU time (seconds) & $err_G$   \\
\hline   Implicit filter ($4,000$ state points )&   $83.14$  & $0.0924$\\
\hline   Particle filter ($15,000$ particles) &  $82.89$ &  $0.2545$\\
\hline   Particle filter ($20,000$ particles) &  $ 142.61$ & $0.1687$\\
\hline   Particle filter ($25,000$ particles) &  $ 209.27$  & $0.1057$\\
\hline
\end{tabular}\end{center}
\end{table}

\section{Conclusions}\label{sec:con}
In this work, we proposed an efficient meshfree implicit filter algorithm by evaluating the conditional PDF on meshfree points in the state space.  These meshfree points are chosen adaptively according to the system state evolution. We also apply Shepard's method as the meshfree interpolation method to compute interplants with random state points. In order to address the degeneracy of the random points, we use importance sampling method to construct a resample step. Numerical examples demonstrate the effectiveness and efficiency of our algorithm. In the future, we plan to perform a rigorous numerical analysis for the meshfree implicit filter algorithm. 

\begin{thebibliography}{10}

\bibitem{Baker2013}
Syed Baker, Hart Poskar, Falk Schreilber, and Bjorn Junker.
\newblock An improved constraint filtering technique for inferring hidden
  states and parameters of a biological model.
\newblock {\em Bioinformatics}, 29:1052--1059, 2013.

\bibitem{Bao-implicit}
F.~Bao, Y.~Cao, and X.~Han.
\newblock An implicit algorithm of solving nonlinear filtering problems.
\newblock {\em Commun. Comput. Phys.}, 16(2):382--402, 2014.

\bibitem{EKF-tracking}
Yaakov Bar-Shalom and Thomas~E. Fortmann.
\newblock {\em Tracking and data association}, volume 179 of {\em Mathematics
  in Science and Engineering}.
\newblock Academic Press Inc., San Diego, CA, 1988.

\bibitem{Bensoussan2009}
Alain Bensoussan, Jussi Keppo, and Suresh~P. Sethi.
\newblock Optimal consumption and portfolio decisions with partially observed
  real prices.
\newblock {\em Math. Finance}, 19(2):215--236, 2009.

\bibitem{particle-filter-resample}
Miodrag Boli{\'c}, Petar~M. Djuri{\'c}, and Sangjin Hong.
\newblock Resampling algorithms and architectures for distributed particle
  filters.
\newblock {\em IEEE Trans. Signal Process.}, 53(7):2442--2450, 2005.

\bibitem{Budhiraja-Survey}
Amarjit Budhiraja, Lingji Chen, and Chihoon Lee.
\newblock A survey of numerical methods for nonlinear filtering problems.
\newblock {\em Phys. D}, 230(1-2):27--36, 2007.

\bibitem{Crisan-PF}
D.~Crisan and O.~Obanubi.
\newblock Particle filters with random resampling times.
\newblock {\em Stochastic Process. Appl.}, 122(4):1332--1368, 2012.

\bibitem{Crisan-Xiong-PF}
Dan Crisan and Jie Xiong.
\newblock A central limit type theorem for a class of particle filters.
\newblock {\em Commun. Stoch. Anal.}, 1(1):103--122, 2007.

\bibitem{cd2002}
Crisan D. and Doucet A.
\newblock A survey of convergence results on particle filtering methods for
  practitiners.
\newblock {\em IEEE Trans. Sig. Proc.}

\bibitem{Dun-EKF}
Jind{\v{r}}ich Dun{\'{\i}}k, Miroslav {\v{S}}imandl, and Ond{\v{r}}ej Straka.
\newblock Unscented {K}alman filter: aspects and adaptive setting of scaling
  parameter.
\newblock {\em IEEE Trans. Automat. Control}, 57(9):2411--2416, 2012.

\bibitem{Elliott2013}
Robert~J. Elliott and Tak~Kuen Siu.
\newblock Option pricing and filtering with hidden {M}arkov-modulated pure-jump
  processes.
\newblock {\em Appl. Math. Finance}, 20(1):1--25, 2013.

\bibitem{Fasshauer2007}
Gregory~F. Fasshauer.
\newblock {\em Meshfree Approximation Methods with MATLAB}.
\newblock Interdisciplinary Mathematical Sciences (Book 6). World Scientific
  Publishing Company, 2007.

\bibitem{particle-filter}
N.J Gordon, D.J Salmond, and A.F.M. Smith.
\newblock Novel approach to nonlinear/non-gaussian bayesian state estimation.
\newblock {\em IEE PROCEEDING-F}, 140(2):107--113, 1993.

\bibitem{Hairer2011}
M.~Hairer, A.~Stuart, and J.~Voss.
\newblock Signal processing problems on function space: {B}ayesian formulation,
  stochastic {PDE}s and effective {MCMC} methods.
\newblock {\em The {O}xford handbook of nonlinear filtering}, pages 833--873,
  2011.

\bibitem{Yang2012}
G.~Huang and P.~Mehta.
\newblock Joint probabilistic data association-feedback particle  
filter with
  applications to multiple target tracking.
\newblock {\em Procs. of American Control Conference}, 2012.

\bibitem{EKF}
Adrew~H. Jazwinski.
\newblock {\em Stochastic Processing and Filtering Theory}, volume~64.
\newblock Academic Press, New York, 1973.

\bibitem{Julier-EKF}
Simon~J. Julier and Joseph~J. LaViola, Jr.
\newblock On {K}alman filtering with nonlinear equality constraints.
\newblock {\em IEEE Trans. Signal Process.}, 55(6, part 2):2774--2784, 2007.

\bibitem{Kim2013}
D.~Kim, T.~Song, and D.~Musicki.
\newblock Highest probability data association for multi-target particle
  filtering with nonlinear measurements.
\newblock {\em IEICE TRANSACTIONS on Communications}, E96-B (1):281--290, 2013.

\bibitem{Kulikov-EKF}
Gennady~Yu. Kulikov and Maria~V. Kulikova.
\newblock Accurate numerical implementation of the continuous-discrete extended
  {K}alman filter.
\newblock {\em IEEE Trans. Automat. Control}, 59(1):273--279, 2014.

\bibitem{Little-Jones2013}
Max~A. Little and Nick~S. Jones.
\newblock Signal processing for molecular and cellular biological physics: an
  emerging field.
\newblock {\em Philos. Trans. R. Soc. Lond. Ser. A Math. Phys. Eng. Sci.},
  371(1984):20110546, 18, 2013.

\bibitem{2D_PopulationModel}
Benzekry Sebastien.
\newblock Mathematical analysis of a two-dimensional population model of
  metastatic growth including angiogenesis.
\newblock {\em J. Evol. Equ.}, 11:187--213, 2011.

\bibitem{Singh2013}
S.~Singh, S.~Digumarthy, A.~Back, J.~Shepard, and Kalra M.
\newblock Radiation dose reduction for chest ct with non-linear adaptive
  filters.
\newblock {\em Acta Radiologica}, 55(2):169--174, 2013.

\bibitem{Stannat2011}
W.~Stannat.
\newblock Stability of the optimal filter for nonergodic signals---a
  variational approach.
\newblock {\em The Oxford handbook of nonlinear filtering}, pages 374--399,
  2011.

\end{thebibliography}

\def\cprime{$'$}

\end{document}